\documentclass[11pt]{article}
\usepackage[matrix,arrow,curve,cmtip]{xy}
\usepackage{amssymb}
\usepackage{theorem}
\usepackage{latexsym}
\usepackage{graphicx}
\usepackage{theorem}
\usepackage[a4paper,width=140mm,height=220mm]{geometry}

\def\to{\mbox{$\xymatrix@1@C=5mm{\ar@{->}[r]&}$}}
\def\tto{\mbox{$\xymatrix@1@C=5mm{\ar@{=>}[r]&}$}}

\def\bkar{\mbox{$\xymatrix@1@C=5mm{\ar@{->}[l]&}$}}
\def\distsign{\begin{picture}(0,0)\put(0,0){\circle{4}}\end{picture}}
\def\dist{\mbox{$\xymatrix@1@C=5mm{\ar@{->}[r]|{\distsign}&}$}}
\def\bkdist{\mbox{$\xymatrix@1@C=5mm{\ar@{->}[l]|{\distsign}&}$}}
\def\biar{\mbox{$\xymatrix@1@C=5mm{\ar@<1.5mm>[r]\ar@<-0.5mm>[r]&}$}}
\def\bidist{\mbox{$\xymatrix@1@C=5mm{\ar@<1.5mm>[r]|{\distsign}\ar@<-0.5mm>[r]|{\distsign}&}$}}
\def\adjar{\mbox{$\xymatrix@1@C=5mm{\ar@<1.5mm>@{<-}[r]\ar@<-0.5mm>[r]&}$}}
\def\adjdist{\mbox{$\xymatrix@1@C=5mm{\ar@<1.5mm>@{<-}[r]|{\distsign}\ar@<-0.5mm>[r]|{\distsign}&}$}}
\def\iso{\mbox{$\xymatrix@1@C=6mm{\ar@{->}[r]^{\sim}&}$}}
\def\doubiso{\mbox{$\xymatrix@1@C=6mm{\ar@{<->}[r]^{\sim}&}$}}
\def\doubar{\mbox{$\xymatrix@1@C=6mm{\ar@{<->}[r]&}$}}
\def\rotleq#1{\mbox{\rotatebox[origin=c]{#1}{$\leq$}}}
\def\rotgeq#1{\mbox{\rotatebox[origin=c]{#1}{$\geq$}}}

\font\atip = xycmat12 at 10pt
\font\btip = xycmbt12 at 10pt
\def\tip#1{
\mbox{
\begin{picture}(0,0)
\put(-725,-25){\atip\char #1}
\put(-725,-25){\btip\char #1}
\end{picture} }}
\newsavebox{\westsouthwesthead}
\savebox{\westsouthwesthead}{%
\tip{119}}
\newcommand{\wswhead}{\usebox{\westsouthwesthead}}
\def\Endoar#1{
\setlength{\unitlength}{0.01pt}
\ifinner
\mbox{
\begin{picture}(300,1200)(1300,0)
\put(1450,680){\mbox{\footnotesize{$#1$}}}
\put(900,770){\oval(900,900)[t]}
\put(900,770){\oval(900,900)[br]}
\put(900,300){\wswhead}
\end{picture}}
\else
\mbox{
\begin{picture}(200,1400)(1600,400)
\put(2100,1300){\mbox{$#1$}}
\put(1300,1300){\oval(1300,1300)[t]}
\put(1300,1300){\oval(1300,1300)[br]}
\put(1300,600){\wswhead}
\end{picture}}
\fi
\setlength{\unitlength}{1pt}}
\def\endoar{
\setlength{\unitlength}{0.01pt}
\ifinner
\mbox{
\begin{picture}(300,1200)(1300,0)
\put(900,770){\oval(900,900)[t]}
\put(900,770){\oval(900,900)[br]}
\put(900,300){\wswhead}
\end{picture}}
\else
\mbox{
\begin{picture}(200,1400)(1600,400)
\put(1300,1300){\oval(1300,1300)[t]}
\put(1300,1300){\oval(1300,1300)[br]}
\put(1300,600){\wswhead}
\end{picture}}
\fi
\setlength{\unitlength}{1pt}}
\def\Endodist#1{
\setlength{\unitlength}{0.01pt}
\ifinner
\mbox{
\begin{picture}(300,1200)(1300,0)
\put(1600,770){\mbox{\footnotesize{$#1$}}}
\put(900,770){\oval(900,900)[t]}
\put(900,770){\oval(900,900)[br]}
\put(1300,970){\circle{400}}
\put(900,300){\wswhead}
\end{picture}}
\else
\mbox{
\begin{picture}(200,1400)(1600,400)
\put(2100,1300){\mbox{$#1$}}
\put(1300,1300){\oval(1300,1300)[t]}
\put(1300,1300){\oval(1300,1300)[br]}
\put(1850,1600){\circle{400}}
\put(1300,600){\wswhead}
\end{picture}}
\fi
\setlength{\unitlength}{1pt}}

\newtheorem{theorem}{Theorem}[section]
\newtheorem{lemma}[theorem]{Lemma}
\newtheorem{definition}[theorem]{Definition} 
\newtheorem{proposition}[theorem]{Proposition}
\newtheorem{corollary}[theorem]{Corollary}

{\theorembodyfont{\upshape}\newtheorem{example}[theorem]{Example}}
\newcommand{\proof}{\noindent {\it Proof\ }: }
\def\endofproof{$\mbox{ }\hfill\Box$\par\vspace{1.8mm}\noindent}

\def\+{^{\dagger}}
\def\etal{~{\it et~al.}}

\def\Cont{{\sf Cont}}
\def\Cocont{{\sf Cocont}}

\def\Rel{{\sf Rel}}
\def\Inf{{\sf Inf}}

\def\:{\colon}
\def\1{{\bf 1}}
\def\impl{\Rightarrow}
\def\2{{\bf 2}}
\def\3{{\bf 3}}

\def\Quant{{\sf Quant}}
\def\Set{{\sf Set}}
\def\QUANT{{\sf QUANT}}

\def\Nat{{\sf Nat}}

\def\Lax{{\sf Lax}}
\def\OpLax{{\sf OpLax}}
\def\skel{_{\sf skel}}
\def\cc{_{\sf cc}}

\def\op{^{\sf op}}

\def\co{^{\sf co}}

\def\dom{{\sf dom}}
\def\cod{{\sf cod}}

\def\CAT{{\sf CAT}}
\def\Ab{{\sf Ab}}
\def\Set{{\sf Set}}

\def\Sup{{\sf Sup}}
\def\Cat{{\sf Cat}}
\def\Matr{{\sf Matr}}

\def\Bim{{\sf Bim}}

\def\Dist{{\sf Dist}}

\def\Cat{{\sf Cat}}
\def\CAT{{\sf CAT}}

\def\Map{{\sf Map}}

\def\Q{{\cal Q}}

\def\D{{\cal D}}
\def\P{{\cal P}}
\def\V{{\cal V}}
\def\W{{\cal W}}

\def\U{{\cal U}}

\def\colim{\mathop{\rm colim}}
\def\lim{\mathop{\rm lim}}
\def\bbA{\mathbb{A}}
\def\bbB{\mathbb{B}}
\def\bbC{\mathbb{C}}
\def\bbD{\mathbb{D}}

\def\bbM{\mathbb{M}}
\def\bbN{\mathbb{N}}
\def\bbF{\mathbb{F}}

\def\bbI{\mathbb{I}}

\def\bbR{\mathbb{R}}

\def\tensor{\otimes}
\def\<{\langle}
\def\>{\rangle}

\title{Categorical structures enriched in a quantaloid: \\ categories, distributors and
functors}
\author{Isar Stubbe\footnote{D\'epartement de Math\'ematique, Universit\'e de Louvain, Chemin du Cyclotron 2,
1348 Louvain-la-Neuve (Belgique), {\tt
i.stubbe@math.ucl.ac.be}.}}
\date{May 10, 2004}

\begin{document}

\maketitle

\begin{abstract} We thoroughly treat several familiar and less familiar
definitions and results concerning categories, functors and distributors enriched in a base quantaloid $\Q$. In analogy with $\V$-category theory we discuss such things as adjoint functors, (pointwise) left Kan extensions, weighted (co)limits, presheaves and free (co)com\-ple\-tion, Cauchy completion and Morita equivalence. With an appendix on the universality of the quantaloid $\Dist(\Q)$ of $\Q$-enriched categories and distributors.
\end{abstract}

\section{Introduction}

The theory of categories enriched in a symmetric monoidal closed category $\V$ is, by now, well known [B\'enabou, 1963, 1965; Eilenberg and Kelly, 1966; Lawvere, 1973; Kelly, 1982]. For such a $\V$ with ``enough'' (co)limits the theory of $\V$-categories, distributors and functors can be pushed as far as needed: it includes such things as (weighted) (co)limits in a $\V$-category, $\V$-presheaves on a $\V$-category, Kan extensions of enriched functors, Morita theory for $\V$-categories, and so on.
\par
Monoidal categories are precisely one-object bicategories [B\'enabou, 1967]. It is thus natural to ask in how far $\V$-category theory can be generalized to $\W$-category theory, for $\W$ a general bicategory. But, whereas in  $\V$-category theory one usually assumes the symmetry of the tensor in $\V$ (which is essential to show that $\V$ is itself a $\V$-category with hom-objects given by the right adjoint to tensoring), in working over a general bicategory $\W$ we will have to sacrifice this symmetry: tensoring objects in $\V$ corresponds to composing morphisms in $\W$ and in general it simply does not make sense for the composition $g\circ f$ of two arrows $f,g$ to be ``symmetric''.
\par
On the other hand, we can successfully translate the notion of closedness of a monoidal category $\V$ to the more general setting of a bicategory $\W$: ask that, for any object $X$ of $\W$ and any arrow $f:A\to B$ in $\W$, both functors 
\begin{eqnarray}
\label{i1} & -\circ f:\W(B,X)\to\W(A,X):x\mapsto x\circ f, & \\
\label{i2} & f\circ -:\W(X,A)\to\W(X,B):x\mapsto f\circ x &
\end{eqnarray}
have respective right adjoints
\begin{eqnarray}
\label{i3} & \{f,-\}:\W(A,X)\to\W(B,X):y\mapsto\{f,y\}, & \\
\label{i4} & [f,-]:\W(X,B)\to\W(X,A):y\mapsto[f,y]. &
\end{eqnarray}
Such a bicategory $\W$ is said to be closed. Some call an arrow such as $\{f,y\}$ a (right) extension and $[f,y]$ a (right) lifting (of $y$ through $f$).
\par
Finally, by saying that $\V$ has ``enough limits and colimits'' is in practice often meant that $\V$ has small limits and small colimits. In a bicategory $\W$ the analogue is straightforward: now ask for $\W$ to have in its hom-categories small limits and small colimits (i.e.~$\W$ is locally complete and cocomplete).
\par
So, to summarize, when trying to develop category theory over a base bicategory $\W$ rather than a base monoidal category $\V$, it seems reasonable to work with a base bicategory which is closed, locally complete and locally cocomplete.
Note that in such a bicategory $\W$, due to its closedness, composition always distributes on both sides over colimits of morphisms:
\begin{eqnarray}
\label{i5} & f\circ(\colim_{i\in I} g_i)\cong \colim_{i\in I}(f\circ g_i), & \\
\label{i6} & (\colim_{j\in J} f_j)\circ g\cong \colim_{j\in J}(f_j\circ g). & 
\end{eqnarray}
That is to say, the local colimits are stable under composition. (But this does not hold in general for local limits!)
\par
We will focus on a special case of these closed, locally complete and locally cocomplete bicategories: namely, we study such bicategories whose hom-categories are moreover small and skeletal. Thus the hom-categories are simply complete lattices. We will write the local structure as an order, and local limits and colimits of morphisms as their infimum, resp.~supremum---so for arrows with same domain and codomain we have things like $f\leq f'$, $\bigvee_{i\in I}f_i$, $\bigwedge_{j\in J}g_j$, etc. In particular (\ref{i5}) and (\ref{i6}) become
\begin{eqnarray}
\label{i5bis} & f\circ(\bigvee_i g_i)= \bigvee_i(f\circ g_i), & \\
\label{i6bis} & (\bigvee_j f_j)\circ g= \bigvee_j(f_j\circ g). & 
\end{eqnarray}
The adjoint functor theorem says that the existence of the adjoints  (\ref{i3}) and (\ref{i4}) to the composition functors (\ref{i1}) and (\ref{i2}) (not only implies but also) is implied by their distributing over suprema of morphisms as in (\ref{i5bis}) and (\ref{i6bis}). Such bicategories -- whose hom-categories are complete lattices and whose composition distributes on both sides over arbitrary suprema -- are called quantaloids. A one-object quantaloid is a quantale\footnote{It was C.~Mulvey [1986] who introduced the word `quantale' in his work on (non-commutative) $C^*$-algebras as a contraction of `quantum' and `locale'.}. So a quantaloid $\Q$ is a $\Sup$-enriched category (and a quantale is monoid in $\Sup$).
\par
We will argue below that ``$\V$-category theory'' can be generalized to ``$\Q$-category theory'', where now $\Q$ denotes a quantaloid. This is a particular case of the theory of ``$\W$-category theory'' as pioneered by [B\'enabou, 1967;   Walters, 1981; Street, 1983a]. But we feel that this particular case is also of particular interest: many examples of bicategory-enriched categories are really quantaloid-enriched. In our further study of categorical structures enriched in a quantaloid we rely heavily on the basic theory of $\Q$-categories, distributors and functors; however, often we could not find an appropriate reference for one or another basic fact. With this text we wish to provide such a reference. 
\par
A first further study, in [Stubbe, 2004a], is concerned with ``variation and enrichement'', in the sense of [Betti\etal, 1983; Gordon and Power, 1997, 1999]. This is a development of the notion of weighted colimit in a $\Q$-enriched category; in particular, a tensored and cotensored $\Q$-category has all weighted colimits (as in section \ref{abc} below) if and only if it has all conical colimits, if and only if its underlying order has suprema of objects of the same type. This allows for a detailed analysis of the biequivalence between $\Q$-modules and cocomplete $\Q$-categories.
\par
The subject of [Stubbe, 2004b] is that of presheaves on $\Q$-semicategories (``categories without units''), along the lines of [Moens\etal, 2002]; it generalizes the theory of regular modules on rings without units. The point is that certain ``good properties'' of the Yoneda embedding for $\Q$-categories (see section \ref{def} below) are no longer valid for $\Q$-semicategories---for example, a presheaf on a $\Q$-semicategory is not canonically the weighted colimit of representables. Enforcing precisely this latter condition defines what is called a ``regular presheaf'' on a $\Q$-semicategory; there is an interesting theory of ``regular $\Q$-semicategories''.
\par
[Borceux and Cruciani, 1998] gives an elementary description of ordered objects in the topos of sheaves on a locale $\Omega$. This turns out to be all about enriched $\Omega$-semicategories that admit an appropriate Cauchy completion. In [Stubbe, 2004c] we describe more generally the theory of $\Q$-semicategories that admit a well-behaved Cauchy completion, that we want to call ``$\Q$-orders'' (based on the material in section \ref{ghi}). Such $\Q$-orders can equivalently be described as categories enriched in the split-idempotent completion of $\Q$, and so provide a ``missing link'' between [Walters, 1981] and [Borceux and Cruciani, 1998].

\section{Quantaloids}\label{1}  

A sup-lattice is an antisymmetrically ordered set\footnote{By an ``ordered set'' we mean a set endowed with a transitive, reflexive relation; it is what is often called a ``preordered set''. We will be explicit when we mean one such relation that is moreover antisymmetric, i.e.~a ``partial order''.}
$(X,\leq)$ for which every subset $X'\subseteq X$ has a supremum
$\bigvee X'$ in $X$. A morphism of sup-lattices $f\:(X,\leq)\to(Y,\leq)$ is a map
$f\:X\to Y$ that preserves arbitrary suprema. It is well-known that sup-lattices and
sup-morphisms constitute a symmetric monoidal closed category
$\Sup$.
\begin{definition}\label{2} A {\em quantaloid} $\Q$ is a $\Sup$-enriched category.
 A {\em homomorphism}
$F\:\Q\to\Q'$ of quantaloids is a
$\Sup$-enriched functor.
\end{definition}
In principle we don't mind a quantaloid having a proper class of objects. Thus quantaloids and homomorphisms form an illegitimate category
$\QUANT$; small quantaloids define a (true) subcategory $\Quant$. A quantaloid with one
object is often thought of as a monoid in $\Sup$, and is called a {\em quantale}.
\par
In elementary terms, a quantaloid $\Q$ is a category whose hom-sets
are actually sup-lattices, in which composition distributes on both sides over arbitrary
suprema of morphisms. In the same vein, a homomorphism
$F\:\Q\to\Q'$ is a functor of (the underlying) categories that preserves arbitrary suprema
of morphisms.
\par For arrows $f\:A\to B$, $g\:B\to C$, $(h_i\:X\to Y)_{i\in I}$ in a quantaloid $\Q$ we
use notations like
$g\circ f\:A\to C$ for composition, and $\bigvee_ih_i\:X\to Y$ and $h_i\leq h_j\:X\biar Y$
for its local structure; $\Q(A,B)$ is the hom-lattice of arrows from $A$ to $B$. The
identity arrow on an object
$A\in\Q$ is written $1_A\:A\to A$. The bottom element of a sup-lattice $\Q(A,B)$ will
typically be denoted by $0_{A,B}$. With these notations for identity and bottom, we can write a
``Kronecker delta''
$$\delta_{A,B}\:A\to B=\left\{\begin{array}{ll} 1_A\:A\to A & \mbox{ if }A=B,\\
0_{A,B}\:A\to B & \mbox{ otherwise.}
\end{array}\right.$$
\par As we may consider $\Sup$ to be a ``simplified version of $\Cat$'', we may regard
quantaloids as ``simplified bicategories''. Notably, a quantaloid
$\Q$ has small hom-categories with stable local colimits, and therefore it is closed.
Considering morphisms
$f\:A\to B$ and $g\:B\to C$ we note the respective adjoints to composition in $\Q$ as $-\circ
f\dashv\{f,-\}$ and
$g\circ -\dashv[g,-]$; that is to say, for any $h\:A\to C$ we have that
$g\circ f\leq h$ iff $\leq [g,h]$ iff $g\leq\{f,h\}$. Note furthermore that every diagram of
2-cells in a quantaloid trivially commutes.
\par
There are a couple of lemmas involving $[-,-]$ and $\{-,-\}$, that hold in any quantaloid, upon which we rely quite often. Let us give a short overview.
\begin{lemma}\label{04} We work in a quantaloid $\Q$. 
\begin{enumerate}
\item For $f\:A\to B$, $g\:B\to C$ and $h\:A\to C$, $[g,h]=\bigvee\{x\:A\to B\mid g\circ x\leq h\}$ and
$\{f,h\}=\bigvee\{y\:B\to C\mid y\circ f\leq h\}$.
\item The following are
equivalent:
\begin{enumerate}
\item $f\dashv g\:B\adjar A$ in $\Q$ (i.e.~$1_A\leq g\circ f$ and $f\circ g\leq 1_B$);
\item $(f\circ -)\dashv (g\circ-)\:\Q(-,B)\adjar\Q(-,A)$ in $\Sup$;
\item $(-\circ g)\dashv (-\circ f)\:\Q(B,-)\adjar\Q(A,-)$ in $\Sup$;
\item $(g\circ -)=[f,-]\:\Q(-,B)\to\Q(-,A)$ in $\Sup$;
\item $(-\circ f)=\{g,-\}\:\Q(B,-)\to\Q(A,-)$ in $\Sup$.
\end{enumerate}
\item An arrow $f\:A\to B$ has a right adjoint if and only if 
$[f,1_B]\circ f=[f,f]$; in this case the right adjoint to $f$ is $[f,1_B]$. Dually, $g\:B\to A$ in $\Q$ has a left adjoint if and only if  
$g\circ \{g,1_B\}=\{g,g\}$; in this case the left adjoint to $g$ is $\{g,1_B\}$.
\item
For $f,f'\:A\biar B$ and $g,g'\:B\biar A$ such that $f\dashv g$ and $f'\dashv g'$, $f\leq f'$ if and only if $g'\leq g$.
\item
Any $f\:A\to B$ induces, for every
$X\in\Q_0$, an adjunction
$$[-,f]\dashv\{-,f\}\:\Q(X,B)\adjar\Q(A,X)\op.$$ By $\Q(A,X)\op$
is meant the sup-lattice $\Q(A,X)$ with opposite order.
\item
For arrows with suitable domain and codomain we have the (dual) identities
$[f,[g,h]]=[g\circ f,h]$ and $\{k,\{m,n\}\}=\{k\circ
m,n\}$, and the (self-dual) identity
$[x,\{y,z\}]=\{z,[x,y]\}$.
\item For arrows with suitable domain and codomain, $[h,g]\circ[g,f]\leq[h,f]$ and dually $\{l,m\}\circ\{k,l\}\leq\{k,m\}$; also
$1_{\dom(f)}\leq[f,f]$ and dually $1_{\cod(f)}\leq\{f,f\}$.
\end{enumerate}
\end{lemma}
\par
$\Sup$ may also be thought of as an ``infinitary version of $\Ab$''; a quantaloid is then
like a ``ring(oid) with infinitary (and idempotent) sum''. This point of view helps to
explain some of the terminology below---especially when we talk about ``matrices with
elements in a quantaloid $\Q$'' in the appendix.
\begin{example}\label{a00}
Any locale $\Omega$ is a (very particular) quantale. Further on we'll denote $\2$ for the two-element boolean algebra. Given a locale $\Omega$, [Walters, 1981] uses the quantaloid of relations in $\Omega$: the objects of $\Rel(\Omega)$ are the elements of $\Omega$, the hom-lattices are given by $\Rel(\Omega)(u,v)=\{w\in\Omega\mid w\leq u\wedge v\}$ with order inherited from $\Omega$, composition is given by infimum in $\Omega$, and the identity on an object $u$ is $u$ itself. This quantaloid is very particular: it equals its opposite.
\end{example}
\begin{example}\label{a01}
The extended non-negative reals $\bbR^+\cup\{+\infty\}$ form a (symmetric) quantale for the opposite order and addition as binary operation [Lawvere, 1973].
\end{example}
\begin{example}\label{a03}
The ideals in a commutative ring $R$ form a (symmetric) quantale.
When $R$ is not commutative, the two-sided ideals still form a (symmetric) quantale. But, as G.~Van den Bossche [1995] points out (but she credits B.~Lawvere for this idea), there is also quite naturally a ``quantaloid of ideals'' containing a lot more information than just the two-sided ideals. Denoting $\Q_R$ for this structure, define that it has two objects, $0$ and $1$; the hom-sup-lattices are
\begin{quote}
$\Q_R(0,0)=$ additive subgroups of $R$ which are ${\cal Z}(R)$-modules,\\
$\Q_R(0,1)=$ left-sided ideals of $R$,\\
$\Q_R(1,0)=$ right-sided ideals of $R$,\\
$\Q_R(1,1)=$ two-sided ideals of $R$,
\end{quote}
with sum of additive subgroups as supremum:
$$\sum_{k\in K}I_k=\{\mbox{finite sums of elements in }\bigcup_{k\in K}I_k\}.$$
Composition in $\Q_R$ is the multiplication of additive subgroups, as in 
$$I\circ J=\{\mbox{finite sums }i_1j_1+...+i_nj_n\mbox{ with all }i_k\in I\mbox{ and all }j_k\in J\};$$
and the identity arrow on $1$ is $R$ and that on $0$ is ${\cal Z}(R)$. (By ${\cal Z}(R)$ we denote the center of $R$). For a commutative $R$, this quantaloid $\Q_R$ is equivalent as $\Sup$-category to the quantale of ideals in $R$.
\end{example}
More examples can be found in the literature, e.g.~[Rosenthal, 1996].

\section{Three basic definitions}\label{2.1}

In the following
$\Q$ always denotes a quantaloid. By a {\em $\Q$-typed set} $X$ we mean a set $X$ to every
element of which is associated an object of $\Q$: for every $x\in X$ there is a $tx$ in $\Q$ (which is called the {\em type} of $x$ in $\Q$). The notation with a ``$t$'' for the types of
elements in a
$\Q$-typed set is generic; i.e.~even for two different
$\Q$-typed sets
$X$ and
$Y$, the type of an $x\in X$ is written
$tx$, and that of a $y\in Y$ is $ty$.
A $\Q$-typed set $X$ is just a way of writing a set-indexed family $(tx)_{x\in X}$ of
objects (i.e.~a small discrete diagram) in $\Q$. If
$\Q$ is a small quantaloid, then a
$\Q$-typed set is an object of the slice category $\Set/\Q_0$. 
\begin{definition}\label{3} A {\em
$\Q$-enriched category} (or {\em $\Q$-category} for short)
$\bbA$ consists of
\begin{itemize}
\item objects: a $\Q$-typed set $\bbA_0$,
\item hom-arrows: for all $a,a'\in\bbA_0$, an arrow
$\bbA(a',a)\:ta\to ta'$ in
$\Q$,
\end{itemize} satisfying
\begin{itemize}
\item composition-inequalities: for all $a,a',a''\in\bbA_0$,
$\bbA(a'',a')\circ\bbA(a',a)\leq\bbA(a'',a)$ in $\Q$,
\item identity-inequalities: for all $a\in\bbA_0$,
$1_{ta}\leq\bbA(a,a)$ in
$\Q$.
\end{itemize}
\end{definition}
\begin{definition} A {\em distributor}
$\Phi\:\bbA\dist\bbB$ between two 
$\Q$-categories is given by
\begin{itemize}
\item distributor-arrows: for all $a\in\bbA_0$, $b\in\bbB_0$, an arrow
$\Phi(b,a)\:ta\to tb$ in $\Q$
\end{itemize} satisfying
\begin{itemize}
\item action-inequalities: for all $a,a'\in\bbA_0$, $b,b'\in\bbB_0$,
$\bbB(b',b)\circ\Phi(b,a)\leq\Phi(b',a)$ and $\Phi(b,a)\circ\bbA(a,a')\leq\Phi(b,a')$ in $\Q$.
\end{itemize}
\end{definition}
\begin{definition}
A {\em functor} $F\:\bbA\to\bbB$ between $\Q$-categories is
\begin{itemize}
\item object-mapping: a map $F\:\bbA_0\to\bbB_0\:a\mapsto Fa$
\end{itemize} satisfying
\begin{itemize}
\item type-equalities: for all $a\in\bbA_0$, $ta=t(Fa)$ in $\Q$,
\item action-inequalities: for all $a,a'\in\bbA_0$,
$\bbA(a',a)\leq\bbB(Fa',Fa)$ in
$\Q$.
\end{itemize}
\end{definition}
\par None of these definitions requires any of the usual diagrammatic axioms (associativity
and identity axioms for the composition in a category, coherence of the action on a
distributor with the composition in its (co)domain category, the functoriality of a functor)
simply because those conditions, which require the commutativity of certain diagrams of
2-cells in the base quantaloid
$\Q$, hold trivially! It is therefore a property of, rather than an extra structure on, a
given set
$\bbA_0$ of elements with types in
$\Q$ together with hom-arrows $\bbA(a',a)\:ta'\to ta$ (for
$a,a'\in\bbA_0$) to be a $\Q$-category; in other words, if these data determine a
$\Q$-category, then they do so in only one way. Similar for distributors and functors: it is
a property of a given collection of arrows
$\Phi(b,a)\:ta\to tb$ whether or not it determines a distributor between
$\Q$-categories $\bbA$ and $\bbB$, as it is a property of an object mapping
$F\:\bbA_0\to\bbB_0$ whether or not it determines a functor.
\par Our $\Q$-categories are by definition small: they have a set of objects. As a
consequence, the collection of distributors between two $\Q$-categories is always a small
set too, and so is the collection of functors between two
$\Q$-categories. (However, we will soon run into size-related trouble: our base quantaloid
$\Q$ having a proper class of objects will conflict with the
$\Q$-categories being small, in particular in matters related to cocompleteness of
$\Q$-categories. When doing categorical algebra over $\Q$ we will therefore suppose that
$\Q$ is small too. But this hypothesis is not necessary at this point.)
\par
Note that for a distributor $\Phi\:\bbA\dist\bbB$, for all $a\in\bbA$ and $b\in\bbB$,
$$\Phi(b,a)=1_{tb}\circ\Phi(b,a)
\leq \bbB(b,b)\circ\Phi(b,a)
\leq\bigvee_{b'\in\bbB}\bbB(b,b')\circ\Phi(b',a)
\leq\Phi(b,a),$$
that is to say, $\bigvee_{b'\in\bbB}\bbB(b,b')\circ\Phi(b',a)
=\Phi(b,a)$; $\bigvee_{a'\in\bbA}\Phi(b,a')\circ\bbA(a',a)=\Phi(b,a)$ is analogous. A $\Q$-category $\bbA$ is itself a distributor from $\bbA$ to $\bbA$; and the identities above become, for all $a,a''\in\bbA$, $\bigvee_{a'\in\bbA}\bbA(a'',a')\circ\bbA(a',a)=\bbA(a'',a)$. These identities allow to say, in \ref{4}, with a suitable definition for the composition of distributors, that ``$\bbA$ is the identity distributor on $\bbA$''.
\par
There is a notational issue that we should comment on. We have chosen to write the composition of arrows in a base quantaloid $\Q$ ``from right to left'': the composite of $f\:X\to Y$ and $g\:Y\to Z$ is $g\circ f\:X\to Z$. Therefore we have chosen to write the hom-arrows in a $\Q$-enriched category $\bbA$ also ``from right to left'': for two objects $a,a'\in\bbA$, the hom-arrow $\bbA(a',a)$ goes from $ta$ to $ta'$. Doing so it is clear that, for example, the composition-inequality in $\bbA$ is written $\bbA(a'',a')\circ\bbA(a',a)\leq\bbA(a'',a)$, with the {\em pivot} $a'$ nicely in the middle, which we find very natural. Our notational conventions are thus basically those of R.~Street's seminal paper [1983a]; other authors have chosen other notations.
\par The next proposition displays the calculus of $\Q$-categories and distributors.
\begin{proposition}\label{4} $\Q$-categories are the objects, and distributors the arrows,
of a quantaloid
$\Dist(\Q)$ in which
\begin{itemize}
\item the composition $\Psi\tensor_{\bbB}\Phi\:\bbA\dist\bbC$ of two distributors
$\Phi\:\bbA\dist\bbB$ and
$\Psi\:\bbB\dist\bbC$ has as distributor-arrows, for
$a\in\bbA_0$ and $c\in \bbC_0$,
$$(\Psi\tensor_{\bbB}\Phi)(c,a)=\bigvee_{b\in\bbB_0}\Psi(c,b)\circ\Phi(b,a);$$
\item the identity distributor on a $\Q$-category $\bbA$ has as distributor-arrows precisely
the hom-arrows of the category $\bbA$ itself, so we simply write it as
$\bbA\:\bbA\dist\bbA$;
\item the supremum $\bigvee_{i\in I}\Phi_i\:\bbA\dist \bbB$ of given distributors
$(\Phi_i\:\bbA\dist\bbB)_{i\in I}$ is calculated elementwise, thus its distributor-arrows
are, for
$a\in\bbA_0$ and $b\in\bbB_0$,
$$(\bigvee_{i\in I}\Phi_i)(b,a)=\bigvee_{i\in I}\Phi_i(b,a).$$
\end{itemize}
\end{proposition}
The proof of the fact that the data above define a quantaloid is straightforward. Actually, $\Dist(\Q)$ is a universal construction on $\Q$ in $\QUANT$: there is a fully faithful homomorphism of quantaloids
\begin{equation}\label{5.0}
\Q\to\Dist(\Q)\:\Big(f\:A\to B\Big)\mapsto\Big((f)\:*_A\dist *_B\Big)
\end{equation}
sending an object $A\in\Q$ to the $\Q$-category with only one object, say $*$, of type $t*=A$, and hom-arrow $1_A$, and a $\Q$-arrow $f\:A\to B$ to the distributor $(f)\:*_A\dist *_B$ whose single element is $f$. This turns out to be the universal direct-sum-and-split-monad completion of $\Q$ in $\QUANT$. The appendix gives details. (Even for a small base quantaloid
$\Q$, $\Dist(\Q)$ has a proper class of objects; so large quantaloids arise as universal constructions on small ones. Therefore we do not wish to exclude large quantaloids {\em a priori}.)
\par
Since $\Dist(\Q)$ is a quantaloid, it is in particular closed; the importance of this fact cannot be overestimated. Let for example 
$\Theta\:\bbA\dist\bbC$ and
$\Psi\:\bbB\dist\bbC$ be distributors between $\Q$-categories, then $[\Psi,\Theta]\:\bbA\dist\bbB$ is the distributor with distributor-arrows, for $a\in\bbA_0$ and $b\in\bbC_0$,
\begin{equation}\label{5.0.1}
\Big[\Psi,\Theta\Big](b,a)=\bigwedge_{c\in
\bbC_0}\Big[\Psi(c,b),\Theta(c,a)\Big],
\end{equation}
where the liftings on the right are calculated in $\Q$. A similar formula holds for $\{-,-\}$.
\par The category of $\Q$-categories and functors is the obvious one.
\begin{proposition}\label{5} $\Q$-categories are the objects, and functors the arrows, of a
category
$\Cat(\Q)$ in which
\begin{itemize}
\item the composition $G\circ F\:\bbA\to\bbC$ of two functors
$F\:\bbA\to\bbB$ and $G\:\bbB\to\bbC$ is determined by the composition of object maps
$G\circ F\:\bbA_0\to\bbC_0\:a\mapsto G(F(a))$;
\item the identity functor $1_{\bbA}\:\bbA\to\bbA$ on a $\Q$-category $\bbA$ is determined
by the identity object map
$1_{\bbA}\:\bbA_0\to\bbA_0\:a\mapsto a$.
\end{itemize}
\end{proposition}
\par  Every functor between $\Q$-categories induces an adjoint pair of distributors, and the
resulting inclusion of the functor category in the distributor category -- although
straightforward -- is a key element for the development of the theory of $\Q$-enriched
categories.
\begin{proposition}\label{6} For $\Q$-categories and functors $F\:\bbA\to\bbB$ and
$G\:\bbC\to\bbB$, the $\Q$-arrows 
$\bbB(Gc,Fa)\:ta\to tc$, one for each $(a,c)\in\bbA_0\times\bbC_0$, determine a distributor\footnote{There is a converse: for any distributor $\Phi\:\bbA\dist\bbC$ there is a -- universal -- way in which $\Phi=\bbB(G-,F-)$ for certain functors $F\:\bbA\to\bbB$ and $G\:\bbC\to\bbB$.}
$\bbB(G-,F-)\:\bbA\dist\bbC$. In particular, for any functor $F\:\bbA\to\bbB$ the
distributors
$\bbB(1_{\bbB}-,F-)\:\bbA\dist\bbB$ and
$\bbB(F-,1_{\bbB}-)\:\bbB\dist\bbA$ are adjoint in the quantaloid $\Dist(\Q)$:
$\bbB(1_{\bbB}-,F-)\dashv\bbB(F-,1_{\bbB}-)$.
\end{proposition}
\proof For any $a,a'\in\bbA_0$ and $c,c'\in\bbC_0$,
$\bbC(c',c)\circ\bbB(Gc,Fa)\circ\bbA(a,a')\leq\bbB(Gc',Gc)\circ\bbB(Gc,Fa)\circ\bbB(Fa,Fa')\leq\bbB(Gc',Fa)$
by functoriality of $F$ and $G$ and composition in $\bbB$. So $\bbB(G-,F-)$ is a distributor from $\bbA$ to $\bbB$.
\par
To see that $\bbB(1_{\bbB}-,F-)\dashv\bbB(F-,1_{\bbB}-)$ in $\Dist(\Q)$, we must check two inequalities: the unit of the adjunction is due to the composition in $\bbB$ and functoriality of $F$,
$\bbB(F-,1_{\bbB}-)\tensor_{\bbB}\bbB(1_{\bbB}-,F-)
=\bbB(F-,F-)\geq\bbA(-,-)$;
the counit follows from the fact that $\{Fa\mid a\in\bbA_0\}\subseteq\bbB_0$ and -- again -- composition in $\bbB$,
$\bbB(1_{\bbB}-,F-)\tensor_{\bbA}\bbB(F-,1_{\bbB}-)
\leq\bbB(1_{\bbB}-,-)\tensor_{\bbB}\bbB(-,1_{\bbB}-)=\bbB(-,-)$.
\endofproof  
In the following we use the abbreviated notations
$\bbB(-,F-)=\bbB(1_{\bbB}-,F-)$ and $\bbB(F-,-)=\bbB(F-,1_{\bbB}-)$.
\begin{proposition}\label{8} Sending a functor to the left adjoint distributor that it
induces, as in
$$\Cat(\Q)\to\Dist(\Q)\:\Big(F\:\bbA\to\bbB\Big)\mapsto
\Big(\bbB(-,F-)\:\bbA\dist\bbB\Big),$$  is functorial. Sending a functor to the right adjoint determines a similar
contravariant functor.
\end{proposition}
\proof
Trivially an identity functor $1_{\bbA}\:\bbA\to\bbA$ is mapped onto the identity distributor $\bbA(-,1_{\bbA}-)=\bbA(-,-)\:\bbA\dist\bbA$. And if $F\:\bbA\to\bbB$ and $G\:\bbB\to\bbC$, then we want
$\bbC(-,G\circ F-)=\bbC(-,G-)\tensor_{\bbB}\bbB(-,F-)$
but also this holds trivially: read the right hand side of the equation as the action of $\bbB$ on $\bbC(-,G-)$.
\endofproof    The category $\Cat(\Q)$ inherits the local structure from the
quantaloid
$\Dist(\Q)$ via the functor $\Cat(\Q)\to\Dist(\Q)$: we put, for two functors
$F,G\:\bbA\biar\bbB$,
$$F\leq G\iff\bbB(-,F-)\leq\bbB(-,G-)\hspace{3mm}\Big(\iff\bbB(G-,-)\leq\bbB(F-,-)\Big).$$ Thus every hom-set
$\Cat(\Q)(\bbA,\bbB)$ is (neither antisymmetrically nor cocompletely) ordered, and
composition in
$\Cat(\Q)$ distributes on both sides over the local order: $\Cat(\Q)$ is a 2-category,
and the functor in
\ref{8} is a 2-functor which is the identity on objects, ``essentially faithful'' (but not
full), and locally fully faithful. (The contravariant version of this functor reverses the local order!)
\par
In general the ``opposite'' of a $\Q$-category $\bbA$ is not again a $\Q$-category, but rather a $\Q\op$-category: of course $\bbA\op$ is defined to have the same $\Q$-typed set of objects as $\bbA$, but the hom-arrows are reversed: for objects $a,a'$ put $\bbA\op(a',a)=\bbA(a,a')$.
Similarly, for a distributor $\Phi\:\bbA\dist\bbB$ between $\Q$-categories we may define an
opposite distributor between the opposite categories over the opposite base---but this distributor will go in the opposite direction: $\Phi\op\:\bbB\op\dist\bbA\op$ is defined by $\Phi\op(a,b)=\Phi(b,a)$. F For $\Phi\leq\Psi\:\bbA\bidist\bbB$ in $\Dist(\Q)$ it is quite obvious that
$\Phi\op\leq\Psi\op\:\bbB\op\bidist\bbA\op$ in $\Dist(\Q\op)$. Finally, for a functor $F\:\bbA\to\bbB$ between $\Q$-categories, the
{\it same} object mapping $a\mapsto Fa$ determines an arrow
$F\op\:\bbA\op\to\bbB\op$ of $\Cat(\Q\op)$. But if $F\leq G\:\bbA\biar \bbB$ in
$\Cat(\Q)$ then $G\op\leq F\op\:\bbA\op\biar\bbB\op$ in
$\Cat(\Q\op)$. It is
obvious that applying the ``{\sf
op}'' twice, always gives back the original structure. 
\begin{proposition}\label{337} ``Taking opposites'' determines isomorphisms 
$$\Dist(\Q)\cong\Dist(\Q\op)\op\mbox{ and }\Cat(\Q)\cong\Cat(\Q\op)\co$$  of 2-categories (where the ``{\sf
co}'' means: reversing order in the homs).
\end{proposition}
These isomorphisms allow us to ``dualize'' all notions and results concerning $\Q$-enriched categories. For example, the dual of `left Kan extension' in $\Cat(\Q)$ is `right Kan extension', and whereas the former help to characterize left adjoints in $\Cat(\Q)$, the latter do the same for right adjoints. As another example, the dual of `weighted colimit' in a $\Q$-category, is `weighted limit'; the former are `preserved' by all left adjoint functors, the latter by right adjoint ones. Some notions are self-dual, like `equivalence' of categories, or `Cauchy complete' category.
\begin{example}\label{a04}
In a quantaloid $\Q$, the arrows whose {\em codomain} is some object $Y$ are the objects of a $\Q$-category that we will denote -- anticipating \ref{102} -- as $\P Y$: put $t(f\:X\to Y)=X$ and $\P Y(f',f)=[f',f]$. Similarly -- and again anticipating further results -- we denote $\P\+ X$ for the $\Q$-category whose objects are $\Q$-arrows with {\em domain} $X$, with types $t(f\:X\to Y)=Y$, and $\P\+ X(f',f)=\{f,f'\}$.
\end{example}
\begin{example}\label{a05}
Recall that $\2$ denotes the 2-element Boolean algebra; $\2$-categories are orders, distributors are ideal relations, and functors are order-preserving maps.
\end{example}
\begin{example}\label{a05.0}
Consider a sheaf $F$ on a locale $\Omega$; it determines a $\Rel(\Omega)$-category $\bbF$ whose objects are the partial sections of $F$, the type of a section $s$ being the largest $u\in\Omega$ on which it is defined, and whose hom-arrows $\bbF(s',s)$ are, for sections $s,s'$ of types $u,u'$, the largest $v\leq u\wedge u'$ on which (restrictions of) $s$ and $s'$ agree [Walters, 1981].
\end{example}
\begin{example}\label{a06}
A category enriched in the quantale $\bbR\cup\{+\infty\}$ (cf.~\ref{a01}) is a ``generalized metric space'' [Lawvere, 1973]: the enrichment itself is a binary distance function taking values in the positive reals. In particular is the composition-inequality in such an enriched category the triangular inequality. A functor between such generalized metric spaces is a distance decreasing application. (These metric spaces are ``generalized'' in that the distance function is not symmetric, that the distance between two points being zero doesn't imply their being identical, and that the distance between two points may be infinite.)
\end{example}
\begin{example}\label{a07}
A (not necessarily commutative) ring $R$ determines a $\Q_R$-enriched category (with $\Q_R$ as in \ref{a03}): denoting it as ${\sf Comm}_R$, its objects of type $0$ are the elements of $R$, its objects of type $1$ are the elements of ${\cal Z}(R)$, and hom-arrows are given by commutators:
${\sf Comm}(r,s)=\{x\in R\mid rx=xs\}$. 
\end{example}

\section{Some direct consequences}

\subsection*{Underlying orders}

For an object $A$ of a quantaloid $\Q$, denote by $*_A$ the one-object $\Q$-category whose hom-arrow is the identity $1_A$. Given a $\Q$-category $\bbA$, the set $\{a\in\bbA_0\mid ta=A\}$ is in bijection with $\Cat(\Q)(*_A,\bbA)$: any such object $a$ determines a ``constant'' functor $\Delta a\:*_A\to\bbA$; and any such functor $F\:*_A\to\bbA$ ``picks out'' an object $a\in\bbA$. We may thus order the
objects of $\bbA$ by ordering the corresponding constant functors; we will speak of the {\em
underlying order}
$(\bbA_0,\leq)$ of the
$\Q$-category
$\bbA$. Explicitly, for two objects $a,a'\in\bbA_0$ we have that $a'\leq a$ if and only if $A:=ta=ta'$ and for all $x\in\bbA_0$, $\bbA(x,a')\leq\bbA(x,a)$ in $\Q$, or equivalently $\bbA(a',x)\geq\bbA(a,x)$, or equivalently $1_A\leq\bbA(a',a)$. Whenever two objects of $\bbA$ are equivalent in
$\bbA\!$'s underlying order ($a\leq a'$ and $a'\leq a$) then we say that they are {\em isomorphic objects} (and write
$a\cong a'$).
\begin{example}
In the $\Q$-category $\P Y$, whose objects are $\Q$-arrows with codomain $Y$ and whose hom-arrows are given by $[-,-]$, the underlying order coincides with the local order in $\Q$. More precisely, for $f,g\:X\biar Y$ it is the same to say that $f\leq g$ in $(\P Y)_0$ as to say that $f\leq g$ in $\Q(X,Y)$.
\end{example}
\par It is immediate that $F\leq G\:\bbA\biar\bbB$ in $\Cat(\Q)$ if and only if, for all $a\in\bbA_0$, $Fa\leq Ga$
in the underlying order of $\bbB$. This says that the local structure in the 2-category
$\Cat(\Q)$ is ``pointwise order''. Equivalently we could have written that $F\leq G$ if and only if 
$1_{ta}\leq \bbB(Fa, Ga)$ for all $a\in\bbA_0$, which exhibits the resemblance with the usual notion of ``enriched natural transformation''.

\subsection*{Adjoints and equivalences}

An arrow $F\:\bbA\to\bbB$ is {\em left adjoint} to an arrow $G\:\bbB\to\bbA$ in
$\Cat(\Q)$ (and $G$ is then {\em right adjoint} to $F$), written $F\dashv G$, if $1_{\bbA}\leq G\circ F$ and $F\circ G\leq 1_{\bbB}$.
Due to the 2-category $\Cat(\Q)$ being locally ordered, we needn't ask any of the usual
triangular coherence diagrams. The unicity of adjoints in the quantaloid $\Dist(\Q)$ and the locally fully faithful $\Cat(\Q)\to\Dist(\Q)$ allow for the following equivalent expression.
\begin{proposition}\label{10.2} $F\:\bbA\to\bbB$ is left adjoint to $G\:\bbB\to\bbA$ in
$\Cat(\Q)$ if and only if $\bbB(F-,-)=\bbA(-,G-)\:\bbB\dist\bbA$ in $\Dist(\Q)$.
\end{proposition}
\par
Further, $F\:\bbA\to\bbB$ is an {\em equivalence} in $\Cat(\Q)$ if there exists a $G\:\bbB\to\bbA$
such that $G\circ F\cong 1_{\bbA}$ and
$F\circ G\cong 1_{\bbB}$ (in which case also $G$ is an equivalence). Again because $\Cat(\Q)$ is locally ordered, this is the same as saying that $F$ is
both left and right adjoint to some $G$. Again the functor $\Cat(\Q)\to\Dist(\Q)$ gives equivalent expressions.
\begin{proposition}\label{10.3}
Given functors $F\:\bbA\to\bbB$ and $G\:\bbB\to\bbA$ constitute an equivalence in $\Cat(\Q)$ if and only if $\bbB(-,F-)\:\bbA\dist\bbB$ and 
$\bbA(-,G)\:\bbB\dist\bbA$ constitute an isomorphism in $\Dist(\Q)$, if and only if $\bbB(F-,-)\:\bbB\dist\bbA$ and $\bbA(-,G-)\:\bbA\dist\bbB$ constitute an isomorphism in $\Dist(\Q)$.
\end{proposition}
\par
For ordinary categories and functors it is well
known that the equivalences are precisely the fully faithful functors which are essentially
surjective on objects. This holds for
$\Q$-enriched categories too: say that a functor $F\:\bbA\to\bbB$ is {\em fully faithful} if
$$\forall\ a,a'\in\bbA_0:\bbA(a',a)=\bbB(Fa',Fa)\mbox{ in }\Q,$$ and that is {\em
essentially surjective on objects} whenever
$$\forall\ b\in\bbB_0,\ \exists\ a\in\bbA_0:Fa\cong b\mbox{ in }\bbB.$$
In fact, $F$ is
fully faithful if and only if the unit of the adjunction
$\bbB(-,F-)\dashv\bbB(F-,-)$ in $\Dist(\Q)$ saturates to an equality; and if $F$ is essentially surjective on objects then necessarily the co-unit of this
adjunction saturates to an equality (but a functor
$F\:\bbA\to\bbB$ for which the co-unit of the induced adjunction is an equality -- which is
sometimes said to have a ``dense image'' -- is not necessarily essentially surjective on
objects).
\begin{proposition}\label{92} An arrow $F\:\bbA\to\bbB$ in $\Cat(\Q)$ is an equivalence if
and only if it is fully faithful and essentially surjective on objects.
\end{proposition}
\proof Suppose that $F$ is an equivalence, with inverse equivalence $G$. Then, by functoriality of $F$ and $G$, and since
$G\circ F\cong 1_{\bbA}$, 
$\bbA(a',a)\leq \bbB(Fa',Fa)\leq\bbA(GFa',GFa)=\bbA(a',a)$ for all $a,a'\in\bbA_0$. Further,
using
$F\circ G\cong 1_{\bbB}$, it follows that $b\cong FGb$, and $Gb\in\bbA_0$ as required.
Conversely,  supposing that $F$ is fully faithful and essentially surjective on objects,
choose -- using the essential surjectivity of $F$ -- for any
$b\in\bbB_0$ one particular object $Gb\in\bbA_0$ for which $FGb\cong b$. Then the object
mapping 
$\bbB_0\to\bbA_0\:b\mapsto Gb$ is fully faithful itself (and therefore also functorial) due
to the fully faithfulness of
$F$: 
$\bbB(b',b)=\bbB(FGb',FGb)=\bbA(Gb',Gb)$ for all $b,b'\in\bbB_0$. It is the required inverse
equivalence.
\endofproof
\par The following is well known for ordinary categories, and holds for $\Q$-categories too; it will be useful further on.
\begin{proposition}\label{482} Suppose that $F\dashv G\:\bbB\adjar\bbA$ in $\Cat(\Q)$. Then
$F$ is fully faithful if and only if
$G\circ F\cong 1_{\bbA}$, and $G$ is fully faithful if and only if $F\circ G\cong 1_{\bbB}$. If moreover $G\dashv H\:\bbA\adjar\bbB$ in $\Cat(\Q)$ then
$F$ is fully faithful if and only if $H$ is fully faithful.
\end{proposition}
\proof
$F\dashv G$ in $\Cat(\Q)$ implies $\bbB(F-,F-)=\bbA(-,G\circ F-)$ in $\Dist(\Q)$; and $F$ is fully faithful if and only if $\bbB(F-,F-)=\bbA(-,-)$ in $\Dist(\Q)$. So, obviously, $F$ is fully faithful if and only if $\bbA(-,G\circ F-)=\bbA(-,-)$ in $\Dist(\Q)$, or equivalently, $G\circ F\cong 1_{\bbA}$ in $\Cat(\Q)$. Likewise for the fully faithfulness of $G$. Suppose now that $F\dashv G\dashv H$, then
$\bbA(G\circ F-,-)=\bbB(F-,H-)=\bbB(-,G\circ H-)$,
which implies that $\bbA(G\circ F-,-)=\bbA(-,-)$ if and only if $\bbB(-,G\circ H-)=\bbA(-,-)$
so the result follows.
\endofproof

\subsection*{Left Kan extensions}

Given $\Q$-categories and functors 
$F\:\bbA\to\bbB$ and $G\:\bbA\to\bbC$, the {\em left Kan extension of $F$ along $G$} is -- in so far it
exists -- a functor
$K\:\bbC\to\bbB$ such that $K\circ G\geq F$ in a universal way: whenever $K'\:\bbC\to\bbB$
satisfies
$K'\circ G\geq F$, then $K'\geq K$.
If the left Kan extension of $F$ along $G$ exists,
then it is essentially unique; we denote it $\<F,G\>$.
So $\<F,G\>\:\bbC\to\bbB$ is the reflection of $F\in\Cat(\Q)(\bbA,\bbB)$ along $$ -\circ G\:\Cat(\Q)(\bbC,\bbB)\to\Cat(\Q)(\bbA,\bbB)\:H\mapsto
H\circ G.$$
\par
Left adjoint functors may be characterized in terms of left Kan extensions; this uses the idea of functors that preserve Kan extensions. Suppose that the left Kan extension $\<F,G\>\:\bbC\to\bbB$ exists; then a functor $F'\:\bbB\to\bbB'$ is said to {\em preserve}
$\<F,G\>$ if $\<F'\circ F,G\>$ exists and is isomorphic to
$F'\circ \<F,G\>$. And $\<F,G\>$ is {\em absolute} if it is preserved by any
$F'\:\bbB\to\bbB'$.
\begin{proposition}\label{127} For a functor
$F\:\bbA\to\bbB$ between $\Q$-categories the following are equivalent:
\begin{enumerate}
\item $F$ has a right adjoint;
\item $\<1_{\bbA},F\>$ exists and is absolute;
\item $\<1_{\bbA},F\>$ exists and is preserved by
$F$ itself.
\end{enumerate} 
In this case, $\<1_{\bbA},F\>$ is the right adjoint of $F$.
\end{proposition}
\proof
For (1$\impl$2), suppose that $F\dashv G$; we'll prove that $G$ is the absolute left Kan extension of $1_{\bbA}$ along $F$. The unit of the adjunction already says that $G\circ F\geq 1_{\bbA}$; and if $K\:\bbB\to\bbA$ is another functor such that $K\circ F\geq 1_{\bbA}$, then necessarily $K=K\circ 1_{\bbB}\geq K\circ(F\circ G)= (K\circ F)\circ G\geq 1_{\bbA}\circ G=G$, using now the unit of the adjunction. So indeed $G\cong\<1_{\bbA},F\>$. Let now $F'\:\bbA\to\bbB'$ be any functor; we claim that $F'\circ G$ is the left Kan extension of $F'\circ 1_{\bbA}$ along $F$. Already $(F'\circ G)\circ F=F'\circ(G\circ F)\geq F'\circ 1_{\bbA}$ is obvious (using the unit of the adjunction); and if $K'\:\bbB\to\bbB'$ is another functor such that $K'\circ F\geq  F'\circ 1_{\bbA}$, then $K'=K'\circ 1_{\bbB}\geq K'\circ (F\circ G)=(K'\circ F)\circ G\geq (F'\circ 1_{\bbA})\circ G=F'\circ G$. Thus, $G$ is the absolute left Kan extension of $1_{\bbA}$ along $F$.\par
(2$\impl$3) being trivial, we now prove (3$\impl$1); it suffices to prove that $F\dashv\<1_{\bbA},F\>$. But the wanted unit $\<1_{\bbA},F\>\circ F\geq 1_{\bbA}$ is part of the universal property of the left Kan extension; and using the hypothesis that $\<1_{\bbA},F\>$ is preserved by $F$ itself, we have $F\circ\<1_{\bbA},F\>=\<F,F\>$ which is smaller than $1_{\bbB}$ by the universal property of $\<F,F\>$ (since $1_{\bbB}\circ F\geq F$).
\endofproof
The dual of this result says that a functor $F\:\bbA\to\bbB$ has a left adjoint in $\Cat(\Q)$ if and only if the {\em right Kan extension} of $1_{\bbA}$ along $F$ exists and is absolute, if and only if this right Kan extension exists and is preserved by $F$. Of course, the definition of {\em right Kan extension} is the dual of that of left Kan extension: given $F\:\bbA\to\bbB$ and $G\:\bbA\to\bbC$ in $\Cat(\Q)$, the right Kan extension of $F$ along $G$ is a functor $(F,G)\:\bbC\to\bbB$ such that $(F,G)\op\:\bbC\op\to\bbB\op$ is the left Kan extension of $F\op$ along $G\op$ in $\Cat(\Q\op)$. In elementary terms: there is a universal inequality $(F,G)\circ G\leq F$ in $\Cat(\Q)$.
\par
After having introduced weighted colimits in a $\Q$-category, we will discuss {\em pointwise} left Kan extensions: particular colimits that enjoy the universal property given above.

\subsection*{Skeletal categories}

To any order corresponds an antisymmetric order by passing to equivalence classes of elements of the
order. For
$\Q$-enriched categories this can be imitated: say that a $\Q$-category $\bbA$ is {\em
skeletal} if no two different objects in $\bbA$ are isomorphic; equivalently this says that, for any $\Q$-category $\bbC$,
$\Cat(\Q)(\bbC,\bbA)$ is an antisymmetric order, i.e.~that 
$$\Cat(\Q)(\bbC,\bbA)\to\Dist(\Q)(\bbC,\bbA)\:F\mapsto\bbA(-,F-)$$ is injective.
Any $\Q$-category is then equivalent to its ``skeletal quotient''.
\begin{proposition}\label{95} For a $\Q$-category $\bbA$, the following data define a
skeletal $\Q$-cate-gory
$\bbA\skel$:
\begin{itemize}
\item[-] objects: $(\bbA\skel)_0=\Big\{ [a]=\{x\in\bbA_0\mid x\cong a\}\Big|
a\in\bbA_0\Big\}$, with type function
$t[a]=ta$;
\item[-] hom-arrows: for any $[a],[a']\in(\bbA\skel)_0$,
$\bbA\skel([a'],[a])=\bbA(a',a)$.
\end{itemize} The object mapping
$[-]\:\bbA\to\bbA\skel\:x\mapsto[x]$ determines an equivalence in
$\Cat(\Q)$.
\end{proposition}
\proof
The construction of $\bbA\skel$ is well-defined, because:
\begin{itemize}
\item as a set the quotient of
$\bbA_0$ by $\cong$ is well-defined ($\cdot\cong\cdot$ is an equivalence);
\item $t[a]$ is well-defined (all elements of the equivalence class
$[a]$ necessarily have the same type since they are isomorphic);
\item for any $a'\cong b'$ and $a\cong b$ in $\bbA$,
$\bbA(b',b)=\bbA(a',a)$, so $\bbA\skel([a'],[a])$ is well-defined.
\end{itemize} 
$\bbA\skel$ inherits ``composition'' and ``identities'' from
$\bbA$, so it is a $\Q$-category. And $\bbA\skel$ is skeletal: if $[a]\cong[a']$ in $(\bbA\skel)_0$ then necessarily $a\cong a'$ in $\bbA_0$, which by definition implies that $[a]=[a']$.
\par 
The mapping 
$\bbA_0\to(\bbA\skel)_0\:x\mapsto[x]$ determines a
functor which is obviously fully faithful and essentially surjective on objects, so by~\ref{92} it is an equivalence in $\Cat(\Q)$.
\endofproof
Skeletal $\Q$-categories can be taken as objects of a full
 sub-2-category of $\Cat(\Q)$, respectively a full sub-quantaloid of
$\Dist(\Q)$; the
proposition above then says that $\Cat\skel(\Q)\to\Cat(\Q)$ is a biequivalence of 2-categories, and
$\Dist\skel(\Q)\to\Dist(\Q)$ is then an equivalence of quantaloids.
$\Cat\skel(\Q)$ is locally antisymmetrically ordered, so the obvious 2-functor
$$\Cat\skel(\Q)\to\Dist\skel(\Q)\:\Big(F\:\bbA\to\bbB\Big)\mapsto
\Big(\bbB(-,F-)\:\bbA\dist\bbB\Big)$$ is the identity on objects, faithful
(but not full), and locally fully faithful.
\begin{example}
The $\Q$-categories $\P Y$ and $\P\+ X$ as defined in \ref{a04} are skeletal.
\end{example}
\begin{example}
Taking $\Q=\2$, the skeletal $\2$-categories are precisely the antisymmetric orders, i.e.~the partial orders (cf.~\ref{a05}).
\end{example}
\begin{example}
The $\Rel(\Omega)$-category associated to a sheaf on $\Omega$, as in \ref{a05.0}, is skeletal.
\end{example}
\par
Later on we will encounter important $\Q$-categories which are always skeletal: the categories $\P\bbA$ and $\P\+\bbA$ of (contravariant and covariant) presheaves on a given $\Q$-category $\bbA$, and also $\bbA$'s Cauchy completion $\bbA\cc$. (The reason, ultimately, that these $\Q$-categories are skeletal, is that the hom-objects of the quantaloid $\Q$ are antisymmetrically ordered.)

\section{Weighted (co)limits}\label{abc}

\subsection*{Colimits}

We consider a functor $F\:\bbA\to\bbB$ and a distributor $\Theta\:\bbC\dist\bbA$ between $\Q$-categories. Since
$F$ determines a distributor $\bbB(F-,-)\:\bbB\dist\bbA$, and by closedness of the quantaloid
$\Dist(\Q)$, we can calculate the universal lifting $[\Theta,\bbB(F-,-)]\:\bbB\dist\bbC$.
A functor
$G\:\bbC\to\bbB$ is the {\em
$\Theta$-weighted colimit of $F$} if it represents the universal lifting:
$\bbB(G-,-)=[\Theta,\bbB(F-,-)]$. If the $\Theta$-weighted colimit of $F$ exists, then it is necessarily
essentially unique. It therefore makes sense to speak of ``the'' colimit and to denote it by
$\colim(\Theta,F)$; its universal property is thus that
$$\bbB\Big(\colim(\Theta,F)-,-\Big)=\Big[\Theta,\bbB(F-,-)\Big]\mbox{ in }\Dist(\Q).$$
The following diagrams picture the situation:
$$\xy
\xymatrix@=15mm{\bbA\ar[r]^F & \bbB \\
\bbC\ar[u]^{\Theta}|{\distsign}}
\endxy
\hspace{10mm}
\xy
\xymatrix@=15mm{\bbA & \bbB\ar[l]|{\distsign}_{\bbB(F-,-)}\ar[dl]|{\distsign}^{\Big[\Theta,\bbB(F-,-)\Big]} \\
\bbC\ar[u]^{\Theta}|{\distsign}}
\afterPOS{\POS*{\rotleq{90}}\restore}<6mm,-8mm>
\endxy
\hspace{10mm}
\xymatrix@=15mm{\bbA\ar[r]^F & \bbB \\
\bbC\ar[u]^{\Theta}|{\distsign}\ar@{.>}[ur]_{\colim(\Theta,F)}}$$
\par  We wish to speak of ``those $\Q$-categories that admit all weighted colimits'',
i.e.~{\em cocomplete}
$\Q$-categories. But there is a small problem (the
word is well-chosen).
\begin{lemma}\label{112} A cocomplete $\Q$-category $\bbB$ has at least as many objects as
the base quantaloid $\Q$.
\end{lemma}
\proof For each
$X\in\Q$, consider the empty diagram in $\bbB$ weighted by the empty distributor with domain
$*_X$. Then $\colim(\emptyset,\emptyset)\:*_X\to\bbB$ must exist by cocompleteness of
$\bbB$; it ``picks out'' an object of type $X$ in $\bbB$---thus
$\bbB$ must have such an object in the first place.
\endofproof 
 So, by the above, would the base quantaloid have a proper class of
objects, then so would all cocomplete
$\Q$-enriched categories.
This is a problem in the framework of this text, because we didn't even bother
defining ``large'' $\Q$-categories, let alone develop a theory of $\Q$-categories that is
sensitive to these size-related issues. {\it Therefore, from now on, we are happy to work with a small base
quantaloid $\Q$}. The problem of ``small''
versus ``large'' then disappears: also ``small'' $\Q$-categories can be cocomplete.
\par  Next up is a collection of lemmas that will help us calculate colimits.
(There  will be some abuse of notation: for a distributor $\Phi\:*_A\dist\bbA$ and a functor $F\:\bbA\to\bbB$, $\colim(\Phi,F)$ is in principle a functor from $*_A$ to $\bbB$. But such a functor simply ``picks out'' an object of type $A$ in $\bbB$. Therefore we will often think of $\colim(\Phi,F)$ just as being that object. Of course, when the domain of the weight has more than one object, then any colimit with that weight is {\em really} a functor!)
\begin{lemma}\label{110}
\begin{enumerate}
\item For $(\Phi_i\:\bbC\dist\bbA)_{i\in I}$ and $F\:\bbA\to\bbB$, if all colimits involved exist, then $\colim(\bigvee_i\Phi_i,F)$ is the supremum of the $\colim(\Phi_i,F)$ in the order $\Cat(\Q)(\bbC,\bbB)$:
$$\colim(\bigvee_i\Phi_i,F)\cong\bigvee_i\colim(\Phi_i,F).$$
\item For any $F\:\bbA\to\bbB$, $\colim(\bbA,F)$ (exists and) is isomorphic to
$F$.
\item For $\Phi\:\bbD\dist\bbC$, $\Theta\:\bbC\dist\bbA$ and $F\:\bbA\to\bbB$, suppose that $\colim(\Theta,F)$ exists; then $\colim(\Phi,\colim(\Theta,F))$ exists if and only if $\colim(\Theta\tensor_{\bbC}\Phi,F)$
exists, in which case they are isomorphic.
\item For $\Theta\:\bbC\dist\bbA$ and $F\:\bbA\to\bbB$,
$\colim(\Theta,F)$ exists if and only if, for all objects $c\in\bbC_0$,
$\colim(\Theta(-,c),F)$ exists; then
$\colim(\Theta,F)(c)\cong\colim(\Theta(-,c),F)$.
\end{enumerate}
\end{lemma}
\proof 
(1) By assumption,
\begin{eqnarray*}
\bbC\Big(\colim(\bigvee_i\Phi_i,F)-,-\Big)
 & = & \Big[\bigvee_i\Phi_i,\bbB(F-,-)\Big] \\
 & = & \bigwedge_i\Big[\Phi_i,\bbB(F-,-)\Big] \\
 & = & \bigwedge_i\bbB\Big(\colim(\Phi_i,F)-,-\Big)
\end{eqnarray*}
from which it follows that $\colim(\bigvee_i\Phi_i,F)$ is the supremum of the $\colim(\Phi_i,F)$ in $\Cat(\Q)(\bbC,\bbB)$.
\par
(2) Trivially, $[\bbA(-,-),\bbB(F-,-)]=\bbB(F-,-)$.
\par
(3) By a simple calculation 
$$\Big[\Phi,\bbB\Big(\colim(\Theta,F)-,-\Big)\Big]
=\Big[\Phi,\Big[\Theta,\bbB(F-,-)\Big]\Big]
=\Big[\Theta\tensor_{\bbC}\Phi,\bbB(F-,-)\Big]$$
so the $\Phi$-weighted colimit of $\colim(\Theta,F)$ and the $\Theta\tensor_{\bbC}\Phi$-weighted colimit of $F$ are the same thing.
\par 
(4) Necessity is easy: $\colim(\Theta,F)\:\bbC\to\bbB$ is a functor satisfying, for all $c\in\bbC_0$ and $b\in\bbB_0$,
$$\bbB\Big(\colim(\Theta,F)(c),b\Big)=\Big[\Theta(-,c),\bbB(F-,b)\Big];$$
this literally says that $\colim(\Theta,F)(c)$ is the $\Theta(-,c)$-weighted colimit of $F$ (which thus exists).
As for sufficiency, we prove that the mapping
$$K\:\bbC_0\to\bbB_0\:c\mapsto\colim(\Theta(-,c),F)$$ is a functor: $\bbC(c',c)\leq\bbB(Kc',Kc)$ because
\begin{eqnarray*}
\bbB(Kc',Kc)
 & = & \bbB(Kc',-)\tensor_{\bbB}\bbB(-,Kc) \\
 & = & \Big\{\bbB(Kc,-),\bbB(Kc',-)\Big\} \\
 & = & \Big\{\Big[\Phi(-,c),\bbB(F-,-)\Big],\Big[\Phi(-,c'),\bbB(F-,-)\Big]\Big\} 
\end{eqnarray*}
and, with slight abuse of notation\footnote{The $\Q$-arrow $\bbC(c',c)\:tc\to tc'$ is thought of as a one-element distributor $*_{tc}\dist *_{tc'}$.},
\begin{eqnarray*}
 & & \bbC(c',c)\leq \Big\{\Big[\Phi(-,c),\bbB(F-,-)\Big],\Big[\Phi(-,c'),\bbB(F-,-)\Big]\Big\} \\
 & & \iff \Phi(-,c')\tensor_{*_{tc'}}\bbC(c',c)\tensor_{*_{tc}}\Big[\Phi(-,c),\bbB(F-,-)\Big]\leq\bbB(F-,-)
\end{eqnarray*}
which holds because $\Phi(-,c')\tensor_{*_tc'}\bbC(c',c)\leq\Phi(-,c)$ and by the universal property of the ``$[-,-]$''. This functor $K$ has, by construction, the required universal property for it to be the $\Phi$-weighted colimit of $F$.
\endofproof 
\begin{lemma}\label{110.0}
For $F\leq G\:\bbA\biar\bbB$ and $\Theta\:\bbC\dist\bbA$, suppose that $\colim(\Theta,F)$ and $\colim(\Theta,G)$ exist; then $\colim(\Theta,F)\leq\colim(\Theta,G)$. In particular, if $F\cong G$ then
$\colim(\Theta,F)$ exists if and only if $\colim(\Theta,G)$ exists, in which case they are isomorphic.
\end{lemma}
\proof
$[-,-]$ is order-preserving in its second variable, so
$$\bbB(\colim(\Theta,F)-,-)=\Big[\Theta,\bbB(F-,-)\Big]\geq\Big[\Theta,\bbB(G-,-)\Big]=\bbB(\colim(\Theta,G)-,-),$$
hence $\colim(\Theta,F)\leq\colim(\Theta,G)$ as claimed. The second part of the statement is similar.
\endofproof 
\par
Let us point out two important consequences.
\begin{corollary}\label{111} A $\Q$-category $\bbB$ is cocomplete if and only if it admits colimits of the identity functor $1_{\bbB}\:\bbB\to\bbB$ weighted by distributors into $\bbB$ whose domains are one-object $\Q$-categories\footnote{Later we will call those distributors {\em contravariant presheaves}, so $\bbB$ is cocomplete if and only if it admits colimits of $1_{\bbB}$ weighted by contravariant presheaves on $\bbB$.}: for every $\phi\:*_B\dist\bbB$ in $\Dist(\Q)$, there exists an object $b\in\bbB$, necessarily of type $tb=B$, such that 
$$\bbB(b,-)=\Big[\phi,\bbB\Big]\mbox{ in $\Dist(\Q)$.}$$
\end{corollary}
\proof
One direction is trivial. For the other, consider a weight $\Theta\:\bbC\to\bbA$ and a functor $F\:\bbA\to\bbB$.
Then
$$\Big[\Theta,\bbB(F-,-)\Big]=\Big[\Theta,\Big[\bbB(-,F-),\bbB\Big]\Big]=
\Big[\bbB(-,F-)\tensor_{\bbA}\Theta,\bbB(1_{\bbB}-,-)\Big]$$
so that $\colim(\Theta,F)\:\bbC\to\bbB$ exists if and only if $\colim(\bbB(-,F-)\tensor_{\bbA}\Theta,1_{\bbB})\:\bbC\to\bbB$ exists. But the latter exists if and only if, for all objects $c\in\bbC$, $$\colim\Big(\bbB(-,F-)\tensor_{\bbA}\Theta(-,c),1_{\bbB}\Big)\:*_{tc}\to\bbB$$ exists---which they do by hypothesis.
\endofproof
\begin{example}
A $\2$-enriched category $\bbA$ is an ordered set, and a distributor $\phi\:*\dist\bbA$ is a down-closed subset of $\bbA$. Since $[\phi,\bbA]$ is the up-closed set of lower bounds of $\phi$, it is representable if and only if the supremum of $\phi$ exists in $\bbA$. Cocompleteness of $\bbA$ as $\2$-category is thus the same thing as cocompleteness of $\bbA$ as order (since it is equivalent to ask for the suprema of all subsets or only the suprema of down-closed sets).
\end{example}
\begin{corollary}\label{1111}
Given $(F_i\:\bbA\to\bbB)_{i\in I}$ in $\Cat(\Q)$, if the $\bigvee_i\bbB(-,F_i-)$-weighted colimit of $1_{\bbB}$ exists, then it is their supremum in the order $\Cat(\Q)(\bbA,\bbB)$:  
$$\bigvee_iF_i\cong\colim\Big(\bigvee_i\bbB(-,F_i-),1_{\bbB}\Big).$$
\end{corollary}
\proof
The colimit, supposed to exist, satisfies
\begin{eqnarray*}
\bbB\Big(\colim(\bigvee_i\bbB(-,F_i-),1_{\bbB})-,-\Big)
 & = & \Big[\bigvee_i\bbB(-,F_i-),\bbB\Big] \\
 & = & \bigwedge_i\Big[\bbB(-,F_i-),\bbB\Big] \\
 & = & \bigwedge_i\bbB(F_i-,-).
\end{eqnarray*}
It follows straightforwardly that $\colim(\bigvee_i\bbB(-,F_i-),1_{\bbB})$ is the supremum of the $F_i$.
\endofproof
This corollary thus says that, if there is a functor $F\:\bbA\to\bbB$ such that $\bbB(F-,-)=\bigwedge_i\bbB(F_i-,-)$ for given $(F_i\:\bbA\to\bbB)_{i\in I}$, then $F$ is the supremum of the $F_i$. In particular, for an object $a\in\bbA$ it follows that $Fa$ is the supremum of the $F_ia$ in the underlying order $(\bbB,\leq)$.
\begin{example}
Concerning the result in \ref{1111}, it is {\em not} true that each supremum in the order $\Cat(\Q)(\bbA,\bbB)$ is necessarily a weighted colimit in $\bbB$. Consider the orders
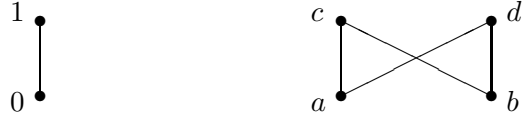
\begin{center}
\setlength{\unitlength}{1 mm}
\begin{picture}(80,30)
\put(10,10){\circle*{1.3}}
\put(10,20){\circle*{1.3}}
\put(10,10){\line(0,1){10}}
\put(50,10){\circle*{1.3}}
\put(50,10){\line(0,1){10}}
\put(50,10){\line(2,1){20}}
\put(70,10){\circle*{1.3}}
\put(50,20){\circle*{1.3}}
\put(70,20){\circle*{1.3}}
\put(70,10){\line(0,1){10}}
\put(70,10){\line(-2,1){20}}
\put(6,8){$0$}
\put(6,20){$1$}
\put(46,8){$a$}
\put(72,8){$b$}
\put(46,20){$c$}
\put(72,20){$d$}
\end{picture}
\end{center}
then the order morphisms $f(0)=a$, $f(1)=c$ and $g(0)=b$, $g(1)=c$ admit the supremum $h(0)=h(1)=c$, but clearly $h(0)$ is not the supremum of $f(0)$ and $g(0)$. This implies that $h$ is not the weighted colimit determined by $f$ and $g$ as in \ref{1111}. 
\end{example}

\subsection*{Limits}

A {\em limit} in a $\Q$-category $\bbA$ is a colimit in the opposite category of $\bbA$. This simple definition
hides the difficulty that $\bbA$ is a $\Q$-category but $\bbA\op$ is a $\Q\op$-category. Let us therefore make the definition of limit in $\bbA$ elementary in terms of $\Q$-enriched structures: for a distributor $\Phi\:\bbB\dist\bbC$ in $\Dist(\Q)$ and a functor $F\:\bbB\to\bbA$ in $\Cat(\Q)$, the $\Phi$-weighted limit of $F$ is -- whenever it exists -- a functor $\lim(\Phi,F)\:\bbC\to\bbA$ in $\Cat(\Q)$ such that
$$\bbA\Big(-,\lim(\Phi,F)-\Big)=\Big\{\Phi,\bbA(-,F-)\Big\}\mbox{ in }\Dist(\Q).$$
These diagrams picture the situation:
$$
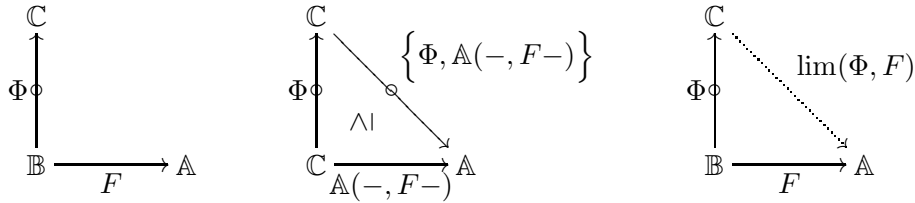
\xy
\xymatrix@=15mm{\bbC\\
\bbB\ar[u]^{\Phi}|{\distsign}\ar[r]_F & \bbA }
\endxy
\hspace{10mm}
\xy
\xymatrix@=15mm{\bbC\ar[dr]|{\distsign}^{\Big\{\Phi,\bbA(-,F-)\Big\}} \\
\bbC\ar[u]^{\Phi}|{\distsign}\ar[r]_{\bbA(-,F-)} & \bbA }
\afterPOS{\POS*{\rotgeq{90}}\restore}<6mm,-14mm>
\endxy
\hspace{10mm}
\xy
\xymatrix@=15mm{\bbC\ar@{.>}[dr]^{\lim(\Phi,F)}\\
\bbB\ar[u]^{\Phi}|{\distsign}\ar[r]_F & \bbA }
\endxy$$
\par
A $\Q$-category $\bbA$ is {\em complete} when the $\Q\op$-category $\bbA\op$ is cocomplete; by \ref{111} we know that this is so precisely when for all $\phi\:*_A\dist\bbA\op$ in $\Dist(\Q\op)$, the $\phi$-weighted colimit of $1_{\bbA\op}$ exists in $\Cat(\Q\op)$; in terms of $\Q$-enriched structures, this occurs precisely when for every\footnote{Later we will call such a distributor a {\em covariant presheaf} on $\bbA$.} $\phi\:\bbA\dist *_A$ in $\Dist(\Q)$, there is an object $a\in\bbA$, necessarily of type $ta=A$, such that
$$\bbA(-,a)=\Big\{\phi,\bbA\Big\}\mbox{ in $\Dist(Q)$.}$$
\begin{proposition}\label{1004}
For any $\Q$-category $\bbA$ and any $\phi\:\bbA\dist *_A$ in $\Dist(\Q)$, $\lim(\phi,1_{\bbA})$ exists if and only if $\colim(\{\phi,\bbA\},1_{\bbA})$ exists, in which case they are isomorphic.
\end{proposition}
\proof
First suppose that $a=\lim(\phi,1_{\bbA})$ exists: $\bbA(-,a)=\{\phi,\bbA\}$. But then, using that $\bbA(-,a)\dashv\bbA(a,-)$ in $\Dist(\Q)$, $\bbA(a,-)=[\bbA(-,a),\bbA]=[\{\phi,\bbA\},\bbA]$, which says that $a$ is the $\{\phi,\bbA\}$-weighted colimit of $1_{\bbA}$.
\par
Conversely, to say that $a=\colim(\{\phi,\bbA\},1_{\bbA})$ exists, means that $\bbA(a,-)=[\{\phi,\bbA\},\bbA]$. Again using the adjunction $\bbA(-,a)\dashv\bbA(a,-)$, and by calculation with $[-,-]$ and $\{-,-\}$ in $\Dist(\Q)$ it follows that $\bbA(-,a)=\{[\{\phi,\bbA\},\bbA],\bbA\}=\{\phi,\bbA\}$, so that indeed $a$ is the $\phi$-weighted colimit of $1_{\bbA}$.
\endofproof
This proposition is valid for categories enriched in any quantaloid $\Q$, hence also for $\Q\op$-categories; translating the result in \ref{1004} for $\Q\op$-categories in terms of $\Q$-categories, we obtain its dual: for any $\psi\:*_A\dist\bbA$ in $\Dist(\Q)$, $\colim(\psi,1_{\bbA})$ exists if and only if $\lim([\psi,\bbA],1_{\bbA})$ exists, in which case they are isomorphic.
\begin{example}
Applied to $\2$-enriched categories, i.e.~ordered sets, \ref{1004} says that the infimum of an up-closed subset equals the supremum of the down-set of its lower bounds (and dually, the supremum of a down-closed subset equals the infimum of the up-closed set of its upper bounds).
\end{example}
The following is then an immediate consequence of \ref{111} and \ref{1004}.
\begin{proposition}\label{1005}
A $\Q$-category is complete if and only if it is cocomplete.
\end{proposition}
This does not mean that ``the theory of complete $\Q$-categories'' is {\em the same} as ``the theory of cocomplete $\Q$-categories'' (although it is its dual): even though the objects are the same, the appropriate morphisms are not. For example, the free cocompletion of a $\Q$-category is different from its free completion (see further).

\subsection*{Absolute (co)limits and (co)continuous functors}

Given a distributor $\Theta\:\bbA\dist\bbB$ and a functor $F\:\bbB\to\bbC$ between $\Q$-categories for which
$\colim(\Theta,F)\:\bbC\to\bbB$ exists, a functor $F'\:\bbB\to\bbB'$ is said to {\em
preserve} this colimit if
$\colim(\Theta, F'\circ F)$ exists and is isomorphic to $F'\circ\colim(\Theta,F)$.
The functor $F'$ is {\em
cocontinuous} if it preserves
all colimits that exist in $\bbB$. And
$\colim(\Theta,F)$ is {\em absolute} if
it is preserved by any functor
$F'\:\bbB\to\bbB'$.
\par
{\em Continuous functors} and {\em absolute limits} are defined dually. For all of the following results we don't bother explicitly writing the dual statements, even though we may use them further on.
\begin{proposition}\label{114} If $\Theta\:\bbC\dist\bbA$ is a left adjoint in
$\Dist(\Q)$, then any colimit with weight $\Theta$ is absolute\footnote{A converse can be
proved too: if any
$\Theta$-weighted colimit is absolute, then $\Theta$ is left adjoint [Street, 1983b].}.
\end{proposition}
\proof If $\Theta\dashv\Psi$ in $\Dist(\Q)$ and $\colim(\Theta,F)$ exists, then necessarily
$$
\bbB\Big(\colim(\Theta,F)-,-\Big)= \Big[\Theta,\bbB(F-,-)\Big]=\Psi\tensor_{\bbA}\bbB(F-,-).$$
From this it is obvious that, for any
$F'\:\bbB\to\bbB'$,
\begin{eqnarray*}
\bbB\Big(F'\circ \colim(\Theta,F)-,-\Big)
 & = & \bbB\Big(\colim(\Theta,F)-,-\Big)\tensor_{\bbB}\bbB'(F-,-)\\
 & = & \Psi\tensor_{\bbA}\bbB(F-,-)\tensor_{\bbB}\bbB'(F'-,-)\\
 & = & \Psi\tensor_{\bbA}\bbB(F'\circ F-,-)\\
 & = & \Big[\Theta,\bbB(F'\circ F-,-)\Big]
\end{eqnarray*} so $\colim(\Theta,F'\circ F)$ exists and is isomorphic to
$F'\circ\colim(\Theta,F)$.
\endofproof
\begin{proposition}\label{116} If $F'\:\bbB\to\bbB'$ is a left adjoint in
$\Cat(\Q)$, then it is cocontinuous\footnote{If $\bbB$ is cocomplete, then also the converse holds---see \ref{132.0}.}.
\end{proposition}
\proof Suppose that $F'\dashv G$ in $\Cat(\Q)$---i.e.~$\bbB'(F'-,-)\dashv\bbB(G-,-)$ in
$\Dist(\Q)$. For any colimit
$$\xymatrix@=15mm{
\bbA\ar[r]^F & \bbB\ar[r]^{F'} & \bbB' \\
\bbC\ar[u]|{\distsign}^{\Theta}\ar@{.>}[ur]_{\colim(\Theta,F)}}$$
we can then calculate
\begin{eqnarray*}
\bbB'\Big(F'\circ\colim(\Theta,F)-,-\Big) 
 & = & \bbB\Big(\colim(\Theta,F)-,-\Big)\tensor_{\bbB}\bbB'(F'-,-) \\
 & = & \Big[\Theta,\bbB(F-,-)\Big]\tensor_{\bbB}\bbB'(F'-,-) \\
 & = & \Big\{\bbB(G-,-),\Big[\Theta,\bbB(F-,-)\Big]\Big\} \\
 & = & \Big[\Theta,\Big\{\bbB(G-,-),\bbB(F-,-)\Big\}\Big] \\
 & = & \Big[\Theta,\bbB(F-,-)\tensor_{\bbB}\bbB'(F'-,-)\Big] \\
 & = & \Big[\Theta,\bbB'(F'\circ F-,-)\Big]
\end{eqnarray*} which shows that
$\colim(\Theta,F'\circ F)$ exists and equals
$F'\circ \colim(\Theta,F)$.
\endofproof
\begin{corollary}\label{116.1} Consider a fully faithful right adjoint $G\:\bbB\to\bbA$ in $\Cat(\Q)$; if $\bbB$ is cocomplete, then so is $\bbA$.
\end{corollary}
\proof
Let us write $F\:\bbA\to\bbB$ for $G$'s left adjoint; now consider any colimit diagram
$$
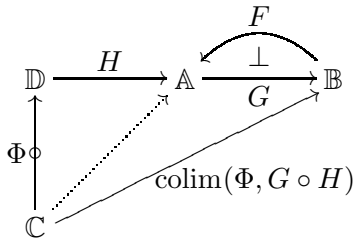
\xymatrix@=15mm{\bbD\ar[r]^H & \bbA\ar[r]_{G}\ar@<3mm>@{}[r]|{\perp} & \bbB\ar@/_6mm/[l]_{F} \\
\bbC\ar[u]|{\distsign}^{\Phi}\ar@<-0.8mm>[urr]_{\ \ \colim(\Phi,G\circ H)}\ar@{.>}[ur]}$$
Then by cocontinuity of $F$ we have $F\circ\colim(\Phi,G\circ H)\cong\colim(\Phi,F\circ G\circ H)$ (thus in particular the existence of the right hand side, the dotted arrow in the diagram), and by fully faithfulness of $G$ (see \ref{482}), $\colim(\Phi,F\circ G\circ H)\cong\colim(\Phi,H)$.
\endofproof
The following two corollaries are now immediate.
\begin{corollary}\label{116.3}
Consider a string of three adjoints in $\Cat(\Q)$,
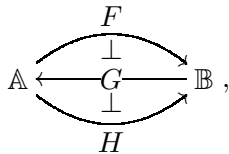
$$\xymatrix@=20mm{
\bbA\ar@/^6mm/[r]^F\ar@/_6mm/[r]_H\ar@<4mm>@{}[r]|{\perp}\ar@<-3mm>@{}[r]|{\perp} & \bbB\ar[l]|{G}},$$
with $F$ fully faithful. If $\bbB$ is cocomplete, then so is $\bbA$.
\end{corollary}
\begin{corollary}\label{116.2} If $F\:\bbA\to\bbB$ is an equivalence in $\Cat(\Q)$, then both $F$ and its inverse equivalence are cocontinuous, and
$\bbA$ is cocomplete if and only if $\bbB$ is.
\end{corollary}

\section{Presheaves}\label{def}

\subsection*{Contravariant presheaves}

Principally the ``contravariant presheaves on a
$\Q$-category
$\bbA$'' should be the objects of a 
$\Q$-category
$\P\bbA$ which ``classifies'' distributors into $\bbA$: 
\begin{equation}\label{1000}
\mbox{for every $\Q$-category
$\bbC$, }\Dist(\Q)(\bbC,\bbA)\simeq\Cat(\Q)(\bbC,\P\bbA).
\end{equation}
(Here we ask {\it a priori} for an equivalence of orders, but it will prove to be an isomorphism of sup-lattices.) Putting $\bbC=*_{C}$ in the above equivalence dictates that a {\em contravariant presheaf of type $C\in\Q$ on $\bbA$} is a distributor $*_C\dist\bbA$.
\par
In what follows we will systematically use small greek letters -- instead of capital ones -- to denote such contravariant presheaves; the elements of such a $\phi\:*_C\dist\bbA$ will be written as
$\phi(a)\:C\to ta$ (instead of $\phi(a,*)$, thus stressing the fact that the action of $*_C$ on $\phi$ is trivial).
\par
The precise sense in which a distributor $\phi\:*_C\dist\bbA$ is {\em contravariant}, is the following: the action of $\bbA$ on $\phi$ can be written as $\bbA(a',a)\leq\{\phi(a),\phi(a')\}$ (for all $a,a'\in\bbA_0$). 
A contravariant presheaf on
$\bbA$ may thus be thought of as a ``functor on $\bbA\op$ with values in
$\Q$'' sending $a$ to $\phi(a)$. To make this precise: indeed, $a\mapsto\phi(a)$ is a functor between $\bbA\op$ and $(\P\+ C)\op$ (see \ref{a04}) in $\Cat(\Q\op)$, or equivalently a functor between $\bbA$ and $\P\+ C$ in $\Cat(\Q)\co$ (by the duality in \ref{337}).
\par
Since we agreed to work over a {\it small} base quantaloid $\Q$, the
collection of contravariant presheaves on a $\Q$-category $\bbA$ forms a
$\Q$-typed set\footnote{Would the objects of the base category $\Q$ form a class, then the
collection of contravariant presheaves on $\bbA$ would be large! The reason -- ultimately --
is that $\P\bbA$ is a cocomplete category---see~\ref{112} and~\ref{107}.}. 
We now need to
enrich this set in $\Q$ in such a way that the classifying property in~(\ref{1000}) holds.
\begin{proposition}\label{102} For any $\Q$-category
$\bbA$, the following defines the $\Q$-category
$\P\bbA$ of contravariant presheaves on
$\bbA$:
\begin{itemize}
\item[-] objects: $(\P\bbA)_0=\{\phi\:*_C\dist\bbA\mbox{ in }\Dist(\Q)\mid C\in\Q\}$ with
types $t(\phi\:*_C\dist \bbA)=C$;
\item[-] hom-arrows: for $\phi\:*_C\dist\bbA$ and $\psi\:*_{C'}\dist\bbA$ put
$\P\bbA(\psi,\phi)\:C\to C'$ to be the single element of the distributor
$[\psi,\phi]\:*_C\dist *_{C'}$.
\end{itemize} 
This $\Q$-category is skeletal, and it satisfies the equivalence in~(\ref{1000}) as follows:
\begin{itemize}
\item for a distributor $\Phi\:\bbC\dist\bbA$ the corresponding functor $F_{\Phi}\:\bbC\to\P\bbA$ maps $c\in\bbC_0$ onto $\Phi(-,c)\:tc\dist\bbA$;
\item and for a functor $F\:\bbC\to\P\bbA$ the corresponding distributor $\Phi_F\:\bbC\dist\bbA$ has elements, for $(c,a)\in\bbC_0\times\bbA_0$, $\Phi_F(a,c)=F(c)(a)$.
\end{itemize}
\end{proposition}
\proof With the properties of ``$[-,-]$'' in $\Dist(\Q)$ it is easily verified that
$\P\bbA$ is a $\Q$-category. For $\phi,\psi\in\P\bbA$, observe that $\phi\leq\psi$ in $(\P\bbA)_0$ if and only if $t\phi=t\psi=:C$ and $(1_C)\leq\P\bbA(\phi,\psi)$, that is to say, $\phi,\psi\:*_C\bidist\bbA$ and $1_C\leq[\phi,\psi]$ in $\Dist(\Q)$, thus precisely when $\phi\leq\psi\:*_C\bidist\bbA$. This shows that the underlying order on the objects of $\P\bbA$ (which are certain distributors) coincides with their order as arrows in $\Dist(\Q)$---so it is antisymmetric, i.e.~$\P\bbA$ is a skeletal $\Q$-category.
As for the equivalence of orders in
(\ref{1000}), it is easily seen that the mappings $F\mapsto \Phi_F$ and $\Phi\mapsto F_{\Phi}$ as in the statement of the proposition define a bijection:
\begin{itemize}
\item for any $c\in\bbC_0$, $F_{\Phi_F}(c)=\Phi_F(-,c)=F(c)$, so $F_{\Phi_F}=F$,
\item for any $(c,a)\in\bbC_0\times\bbA_0$, $\Phi_{F_{\Phi}}(a,c)=F_{\Phi}(c)(a)=\Phi(a,c)$, so $\Phi_{F_{\Phi}}=\Phi$.
\end{itemize}
This bijection preserves order in both directions:
\[\begin{array}{l}
\Phi\leq\Psi\mbox{ in }\Dist(\bbC,\bbA) \\
\iff \forall (c,a)\in\bbC_0\times\bbA_0:\Phi(a,c)\leq\Psi(a,c)\mbox{ in }\Q(tc,ta) \\
\iff \forall c\in\bbC_0: F_{\Phi}(c)\leq F_{\Psi}(c)\mbox{ in }((\P\bbA)_0,\leq) \\
\iff F_{\Phi}\leq F_{\Psi}\mbox{ in }\Cat(\Q)(\bbC,\P\bbA)
\end{array}\]
 and likewise
\[\begin{array}{l}
F\leq G\mbox{ in }\Cat(\Q)(\bbC,\P\bbA) \\
\iff \forall c\in\bbC_0: F(c)\leq G(c)\mbox{ in }((\P\bbA)_0,\leq) \\
\iff \forall (c,a)\in\bbC_0\times\bbA_0:\Phi_F(a,c)\leq\Phi_G(a,c)\mbox{ in }\Q(tc,ta) \\
\iff \Phi_F\leq\Phi_G\mbox{ in }\Dist(\bbC,\bbA).
\end{array}\]
\endofproof
\begin{example}
For a $\Q$-object $Y$, the $\Q$-category of contravariant presheaves on $*_Y$ is denoted $\P Y$ instead of $\P(*_Y)$---see also \ref{a04}.
\end{example}
Very often there will be an abuse of notation: for $\phi\:*_C\dist\bbA$ and $\psi\:*_{C'}\dist\bbA$ there is a {\em distributor} $[\psi,\phi]\:*_C\dist *_{C'}$ with precisely one element (because it goes between one-object categories); therefore we defined $\P\bbA(\psi,\phi)$ to be {\em the single element} of this distributor, for we want $\P\bbA(\psi,\phi)$ to be a $\Q$-arrow. But more often than not we will simply write $\P\bbA(\psi,\phi)=[\psi,\phi]$, thus identifying a quantaloid $\Q$ with the full subcategory of $\Dist(\Q)$ determined by the one-object $\Q$-categories whose hom-arrows are identities (see also (\ref{5.0}) and the appendix).
\par
The equivalence of orders in (\ref{1000}) has a nice side-effect: $\Cat(\Q)(\bbC,\P\bbA)$ is -- by skeletality of $\P\bbA$ -- an antisymmetric order which is order-isomorphic to the sup-lattice $\Dist(\Q)(\bbC,\bbA)$, thus it is a sup-lattice too. In particular are $F\mapsto\Phi_F$ and $\Phi\mapsto F_{\Phi}$ sup-morphisms.
\par
To the identity distributor 
$\bbA\:\bbA\dist\bbA$ corresponds, under the equivalence in~(\ref{1000}), the functor
$$Y_{\bbA}\:\bbA\to\P\bbA:a\mapsto
\Big(\bbA(-,a)\:*_{ta}\dist\bbA\Big)$$
which is, of course, the {\em Yoneda embedding} for
$\bbA$. The presheaf
$\bbA(-,a)\:*_{ta}\dist\bbA$ is said to be {\em represented} by
$a\in\bbA_0$. We can now prove the ``Yoneda lemma'' for presheaves on a $\Q$-enriched
category.
\begin{proposition}\label{104} For any
$\Q$-category
$\bbA$, any object
$a\in\bbA_0$ and any $\phi\:*_C\dist\bbA$ in $\Dist(\Q)$,
$\P\bbA(Y_{\bbA}a,\phi)=\phi(a)$. As a consequence, the Yoneda embedding $Y_{\bbA}\:\bbA\to\P\bbA$ is fully faithful.
\end{proposition}
\proof Using that $\bbA(-,a)\dashv\bbA(a,-)$ in $\Dist(\Q)$, we have by definition of
$\P\bbA(-,-)$ that
$\P\bbA(Y_{\bbA}a,\phi) =[\bbA(-,a),\phi]  =\bbA(a,-)\tensor_{\bbA}\phi  =\phi(a)$. 
\endofproof
\par
An -- almost -- immediate but non the less important consequence of \ref{102} is the following.
\begin{proposition}\label{107} The category
$\P\bbB$ of contravariant presheaves on a $\Q$-category
$\bbB$ is cocomplete.
\end{proposition}
\proof Given a distributor $\Theta\:\bbC\dist\bbA$ and a functor $F\:\bbA\to\P\bbB$, to
$F$ corresponds a unique distributor
$\Phi_F\:\bbA\dist\bbB$ and, in turn, to
$\Phi_F\tensor_{\bbA}\Theta\:\bbC\dist\bbB$ corresponds a unique functor
$F_{\Phi_F\tensor_{\bbA}\Theta}\:\bbC\to\P\bbB$. We claim that this latter functor is the colimit of $F$
weighted by
$\Theta$: for any $c\in\bbC$ and $\beta\in\P\bbB$ we should verify that
$$\P\bbB\Big(F_{\Phi_F\tensor_{\bbA}\Theta}(c),\beta\Big)=\Big[\Phi(-,c),\P\bbB(F-,\beta)\Big].$$
The left hand side equals $[\Phi(-,c),[F-,\beta]]$ and the right hand side equals
$$\Big[(\Phi_F\tensor_{\bbA}\Theta)(-,c),\beta\Big]=\Big[\Phi_F\tensor_{\bbA}\Theta(-,c),\beta\Big]
=\Big[\Theta(-,c),\Big[\Phi_F,\beta\Big]\Big].$$
Since $[F-,\beta]=[\Phi_F,\beta]$, the required universality holds.
\endofproof 
\begin{corollary}\label{109} For any $\phi\:*_C\dist\bbA$, the $\phi$-weighted colimit of $Y_{\bbA}\:\bbA\to\P\bbA$ (exists and) is $\phi$ itself.
\end{corollary}
In other words, every contravariant presheaf on a $\Q$-category $\bbA$ is canonically the colimit of representable presheaves: $\colim(\phi,Y_{\bbA})=\phi$ (and we may write an equality instead of merely an isomorphism, for $\P\bbA$ is skeletal).
\par
Here is another result that should sound familiar.
\begin{proposition}\label{109.1}
The Yoneda embedding $Y_{\bbA}\:\bbA\to\P\bbA$ is continuous.
\end{proposition}
\proof
Consider a weighted limit  (which happens to exist) in $\bbA$, like so:
$$\xymatrix@=15mm{
\bbC\ar@{.>}[dr]^{\lim(\Phi,F)} \\
\bbB\ar[u]|{\distsign}^{\Phi}\ar[r]_F & \bbA\ar[r]_{Y_{\bbA}} & \P\bbA}$$
Then, for any $c\in\bbC_0$,
$$\Big(Y_{\bbA}\circ\lim(\Phi,F)\Big)(c) = \bbA\Big(-,\lim(\Phi,F)(c)\Big)
=\Big\{\Phi(c,-),\bbA(-,F-)\Big\}$$
and therefore
\begin{eqnarray*}
\P\bbA(-,Y_{\bbA}\circ\lim(\Phi,F)-)
 & = & \Big[-,\Big\{\Phi,\bbA(-,F-)\Big\}\Big] \\
 & = & \Big\{\Phi,\Big[-,\bbA(-,F-)\Big]\Big\} \\
 & = & \Big\{\Phi,\P\bbA(-,Y_{\bbA}\circ F-)\Big\}
\end{eqnarray*}
so $Y_{\bbA}\circ \lim(\Phi,F)$ is $\lim(\Phi,Y_{\bbA}\circ F)$.
\endofproof
\par
When $Y_{\bbA}\:\bbA\to\P\bbA$ is a right adjoint, then $\bbA$ is cocomplete (combine \ref{116.1} with \ref{104} and \ref{107}). Also the converse can be shown; actually, $\P\bbA$ is the free cocompletion of $\bbA$ in $\Cat(\Q)$. To make this precise in an elegant way, it is useful to treat pointwise left Kan extensions first.

\subsection*{Pointwise left Kan extensions}

If for $F\:\bbB\to\bbA$ and $G\:\bbB\to\bbC$ in $\Cat(\Q)$ 
the $\bbC(G-,-)$-weighted colimit of $F$ exist, then it is the left Kan extension of $F$ along $G$; so
$\<F,G\>\cong\colim(\bbC(G-,-),F)$ whenever the right hand side exists. 
By lemma~\ref{110} we know that in this case the
action of
$\<F,G\>\:\bbC\to\bbB$ on objects is given by the colimit formula
$\<F,G\>(c)\cong\colim(\bbC(G-,c),F)$. Such a colimit weighted by the right adjoint distributor determined by some functor is therefore said
to be a {\em pointwise left Kan extension}. (Note that, for a pointwise left Kan extension, the two notions of absoluteness, one relevant for Kan extensions and the other relevant for colimits, coincide.)
\par
A $\Q$-category is said to {\em admit all pointwise left Kan extensions} when it admits all
colimits of this special kind.
A cocomplete $\Q$-category then certainly admits all pointwise left Kan extensions.
(For a cocomplete
$\bbB$, any functor
$G\:\bbA\to\bbC$ thus determines an adjunction
$$\<-,G\>\dashv(-\circ G)\:\Cat(\Q)(\bbC,\bbB)\adjar\Cat(\Q)(\bbA,\bbB).$$
This has been called the ``meta-adjointness'' associated
to the concept of Kan extension.) 
But it is also true that a $\Q$-category admitting all pointwise Kan extensions, is cocomplete! This follows from \ref{111} in combination with \ref{122} below; note that the Yoneda embedding plays a crucial r\^ole.
\begin{proposition}\label{122.1}
Suppose that for $F\:\bbB\to\bbA$ and $G\:\bbB\to\bbC$ in $\Cat(\Q)$ the pointwise left Kan extension of $F$ along $G$ exists. If $G$ is fully faithful, then $\<F,G\>\circ G\cong F$.
\end{proposition}
\proof
It is easily seen that, $\<F,G\>$ being a pointwise left Kan extension,
\begin{eqnarray*}
& & \<F,G\>\circ G\cong F \\
& & \iff\bbB(\<F,G\>\circ G-,-)=\bbB(F-,-) \\
& & \iff\bbC(G-,-)\tensor_{\bbC}\bbB(\<F,G\>-,-)=\bbB(F-,-) \\
& & \iff\bbC(G-,-)\tensor_{\bbC}\Big[\bbC(G-,-),\bbB(F-,-)\Big]=\bbB(F-,-) \\
& & \iff\bbC(G-,-)\tensor_{\bbC}\bbC(-,G-)\tensor_{\bbA}\bbB(F-,-)=\bbB(F-,-) \\
& & \iff\bbC(G-,G-)\tensor_{\bbA}\bbB(F-,-)=\bbB(F-,-)
\end{eqnarray*}
so if $G$ is fully faithful, i.e.~$\bbC(G-,G-)=\bbA$, then the above holds.
\endofproof
\begin{proposition}\label{132.0}
Let $\bbA$ be a cocomplete $\Q$-category; then the cocontinuous functors with domain $\bbA$ coincide with the left adjoints.
\end{proposition}
\proof
One direction is a consequence of \ref{116}. To prove the other, 
let $F\:\bbA\to\bbB$ be a cocontinuous functor. By cocompleteness of $\bbA$, $\<1_{\bbA},F\>$ not only exists but is even pointwise, i.e.~it is a colimit, so it must be preserved by $F$. By \ref{127} we conclude that $F$ is a left adjoint.
\endofproof
\begin{proposition}\label{122}
A functor $F\:\bbA\to\bbB$ has a pointwise left Kan extension along $Y_{\bbA}\:\bbA\to\P\bbA$ if and only if, for every $\phi\in\P\bbA$, the $\phi$-weighted colimit of $F$ exists. In this case, $\<F,Y_{\bbA}\>\:\P\bbA\to\bbB$ 
is left adjoint to 
$\<Y_{\bbA},F\>\:\bbB\to\P\bbA$, and $\<F,Y_{\bbA}\>\circ Y_{\bbA}\cong F$. Actually, $\<F,Y_{\bbA}\>$ is then the essentially unique cocontinuous factorization of $F$ through $Y_{\bbA}$.
\end{proposition}
\proof 
If by hypothesis, $\<F,Y_{\bbA}\>=\colim(\P\bbA(Y_{\bbA}-,-),F)$ exists, then 
$\<F,Y_{\bbA}\>(\phi)=\colim(\P\bbA(Y_{\bbA}-,\phi),F)=\colim(\phi,F)$ because $\P\bbA(Y_{\bbA}-,\phi)=\phi$. Conversely, supposing that $\bbB$ allows colimits of $F$ weighted by presheaves on $\bbA$, $\colim(\P\bbA(Y_{\bbA}-,-),F)$ exists because, for all $\phi=\P\bbA(Y_{\bbA}-,\phi)$,
$\colim(\P\bbA(Y_{\bbA}-,\phi),F)$ exists (see~\ref{110}).
\par
Since $\P\bbA$ is cocomplete,
$\<Y_{\bbA},F\>$ exists, is pointwise, and can thus be computed as
$\<Y_{\bbA},F\>(b) \cong \colim(\bbB(F-,b),Y_{\bbA}) \cong \bbB(F-,b)$. It is then a matter of applying~\ref{10.2} to
prove the adjunction:
$$\begin{array}{l}
\<F,Y_{\bbA}\>\dashv\<Y_{\bbA},F\> \\
\iff\forall\phi\in\P\bbA,\forall
b\in\bbB:\bbB\Big(\<F,Y_{\bbA}\>(\phi),b\Big)=\P\bbA\Big(\phi,\<Y_{\bbA},F\>(b)\Big) \\
\iff\forall\phi\in\P\bbA,\forall
b\in\bbB:\bbB\Big(\colim(\phi,F),b\Big)=\Big[\phi,\bbB(F-,b)\Big]
\end{array}$$ which holds by the definition of colimit.
\par
Finally, by fully faithfulness of $Y_{\bbA}$ and \ref{122.1} it follows that $F$ factors through $Y_{\bbA}$ as $\<F,Y_{\bbA}\>$. Would now
$G\:\P\bbA\to\bbB$ be another cocontinuous functor such that
$F\cong G\circ Y_{\bbA}$ then -- by cocontinuity and using the fact that any presheaf
$\phi\in\P\bbA$ is canonically the colimit of representables -- $G$ and $\<F,Y_{\bbA}\>$ agree on all objects of
$\P\bbA$:
\begin{eqnarray*} G(\phi) & = & G(\colim(\phi,Y_{\bbA})) \\
 & \cong & \colim(\phi,G\circ Y_{\bbA}) \\
 & \cong & \colim(\phi, F) \\ 
 & = & \<F,Y_{\bbA}\>(\phi).
\end{eqnarray*}
\endofproof
\begin{corollary}
A $\Q$-category $\bbA$ is cocomplete if and only if $Y_{\bbA}\:\bbA\to\P\bbA$ is a right adjoint in $\Cat(\Q)$.
\end{corollary}
\proof
If $Y_{\bbA}$ is right adjoint, apply \ref{116.1} to $Y_{\bbA}$, using \ref{104} and \ref{107}. If $\bbA$ is cocomplete, apply \ref{122} to the identity functor $1_{\bbA}\:\bbA\to\bbA$: $\<1_{\bbA},Y_{\bbA}\>\dashv\<Y_{\bbA},1_{\bbA}\>=Y_{\bbA}$.
\endofproof

\subsection*{Free cocompletion}

Applying \ref{122} we easily obtain that, for two $\Q$-categories $\bbA$ and $\bbB$, with $\bbB$
cocomplete, any functor
$F\:\bbA\to\bbB$ factors through
$Y_{\bbA}\:\bbA\to\P\bbA$ by the essentially unique cocontinuous -- even left adjoint -- functor
$\<F,Y_{\bbA}\>\:\P\bbA\to\bbB$.
Let now $\Cocont(\Q)$ denote the sub-2-category
of $\Cat(\Q)$ whose objects are the cocomplete $\Q$-categories, and whose morphisms are the
cocontinuous functors, or equivalently (by \ref{132.0}) the left adjoint ones;
there is a forgetful 2-functor
$\U\:\Cocont(\Q)\to\Cat(\Q)$. The above then says that, for $\Q$-categories $\bbA$ and $\bbB$, with $\bbB$ cocomplete,
the Yoneda embedding
$Y_{\bbA}\:\bbA\to\P\bbA$ induces a natural equivalence of orders
$$-\circ Y_{\bbA}\:\Cocont(\Q)(\P\bbA,\bbB)\iso
\Cat(\Q)(\bbA,\bbB)$$ with inverse $F\mapsto\<F,Y_{\bbA}\>$.
\begin{proposition}\label{129}
The presheaf construction provides a left (bi)adjoint to the inclusion of
$\Cocont(\Q)$ in
$\Cat(\Q)$:
$\P\dashv\U\:\Cocont(\Q)\adjar\Cat(\Q)$.
\end{proposition}
\par
For any two $\Q$-categories $\bbA$ and $\bbB$, we have  -- by \ref{102} and the above -- equivalences
\begin{equation}\label{1008}
\Dist(\Q)(\bbA,\bbB)\iso\Cat(\Q)(\bbA,\P\bbB)\iso\Cocont(\Q)(\P\bbA,\P\bbB)
\end{equation}
sending a $\Phi\:\bbA\dist\bbB$ to $\<F_{\Phi},Y_{\bbA}\>\:\P\bbB\to\P\bbA$ (with
$F_{\Phi}\:\bbB\to\P\bbA\:b\mapsto\Phi(-,b)$ as in \ref{102}).
We can prove more.
\begin{proposition}\label{132} The locally ordered\footnote{... but not antisymmetrically so!} category $\Cocont(\Q)$ has stable local colimits, and the action
$$\Dist(\Q)\to\Cocont(\Q)\:\Big(\Phi\:\bbA\dist\bbB\Big)\mapsto\Big(\<F_{\Phi},Y_{\bbA}\>\:\P\bbA\to\P\bbB\Big)$$ is a 2-functor which is locally an equivalence.
\end{proposition}
\proof 
We already know that $\Cocont(\Q)$ is a 2-category; we must show that it has stable local colimits. By \ref{1111} we know how to compute suprema in $\Cat(\Q)(\bbA,\bbB)$ whenever $\bbB$ is cocomplete; we may thus apply this to 
given arrows $(F_i\:\bbA\to\bbB)_{i\in I}$ in $\Cocont(\Q)$:
$$\bigvee_{i\in I}F_i:=\colim\Big(\bigvee_{i\in I}\bbB(-,F_i-),1_{\bbB}\Big).$$
But we must show that $\bigvee_iF_i$ is still a cocontinuous functor (for {\em a priori} it is merely a functor); equivalently (by cocompleteness of $\bbA$), we may assume that each $F_i$ is a left adjoint in $\Cat(\Q)$, and show that $\bigvee_iF_i$ itself is a left adjoint in $\Cat(\Q)$. Let us write $F_i\dashv G_i$ for such adjunctions; then, by \ref{1111} and its dual\footnote{The dual of \ref{1111} says the following: if $\bbA$ is a complete category, then the infimum of functors $(G_i\:\bbB\to\bbA)_{i\in I}$ exists in $\Cat(\Q)(\bbB,\bbA)$, and is given by the $\bigvee_i\bbA(G_i-,-)$-weighted limit of $1_{\bbA}$. This means exactly that this infimum $\bigwedge_iG_i$ satisfies $\bbA(-,\bigwedge_iG_i-)=\bigwedge_i\bbA(-,G_i-)$. By \ref{1005} any complete category is also cocomplete, so this applies in this case!},
$$\bbB(\bigvee_iF_i-,-)=\bigwedge_i\bbB(F_i-,-)=\bigwedge_i\bbA(-,G_i-)=\bbA(-,\bigwedge_iG_i-).$$
That is to say, $\bigvee_iF_i\dashv\bigwedge_iG_i$ in $\Cat(\Q)$. As for the stability of these local colimits, take for example $(F_i\:\bbA\to\bbB)_{i\in I}$ and $G\:\bbB\to\bbC$ in $\Cocont(\Q)$; then
\begin{eqnarray*}
G\circ(\bigvee_iF_i) & = & G\circ\colim\Big(\bigvee_i\bbB(-,F_i-),1_{\bbB}\Big) \\
 & \cong & \colim\Big(\bigvee_i\bbB(-,F_i-),G\Big) \\
 & \cong & \colim\Big(\bbC(-,G-)\tensor_{\bbB}\bigvee_i(-,F_i-),1_{\bbC}\Big) \\
 & \cong & \colim\Big(\bigvee_i(-,G\circ F_i-),1_{\bbC}\Big) \\
 & = & \bigvee_i(G\circ F_i)
\end{eqnarray*}
because $G$ is cocontinuous. Now consider $F\:\bbA\to\bbB$ and $(G_i\:\bbB\to\bbC)_{i\in I}$, then for any $a\in\bbA_0$
\begin{eqnarray*}
\Big((\bigvee_iG_i)\circ F\Big)(a)
 & = & \colim\Big(\bigvee_i\bbC(-,G_i-),1_{\bbC}\Big)(Fa) \\
 & \cong & \colim\Big(\bigvee_i\bbC(-,G_i-)\tensor_{\bbB}\bbB(-,Fa),1_{\bbC}\Big) \\
 & \cong & \colim\Big(\bigvee_i\bbC(-,G_i(Fa)),1_{\bbC}) \\
 & \cong & \colim\Big(\bigvee_i\bbC(-,G_i\circ F-),1_{\bbC})(a) \\
 & = & \Big(\bigvee_i(G_i\circ F)\Big)(a)
\end{eqnarray*}
by properties of the colimit -- see \ref{110} -- so $(\bigvee_iG_i)\circ F=\bigvee_i(G_i\circ F)$.
\par
To prove the functoriality of the action described in the statement of the proposition, first observe that for an identity distributor
$\bbA\:\bbA\dist\bbA$, both $\<Y_{\bbA},Y_{\bbA}\>$ and $1_{\P\bbA}$ are cocontinuous (left adjoint) factorizations of $Y_{\bbA}$ through $Y_{\bbA}$, so -- applying~\ref{129} -- at least they are isomorphic functors; $\P\bbA$ being skeletal implies that this isomorphism is an equality. For $\Phi\:\bbA\dist\bbB$ and $\Psi\:\bbB\dist\bbC$, to check that
$\<F_{\Psi},Y_{\bbB}\>\circ\<F_{\Phi},Y_{\bbA}\> =\<F_{\Psi\tensor_{\bbB}\Phi},Y_{\bbA}\>$ in $\Cat(\Q)(\P\bbA,\P\bbC)$,
it suffices to see -- again by \ref{129} and by skeletality of $\P\bbC$ -- that
$\<F_{\Psi},Y_{\bbB}\>\circ\<F_{\Phi},Y_{\bbA}\>$ is a cocontinuous factorization of $F_{\Psi\tensor_{\bbB}\Phi}$ through $Y_{\bbA}$. As composition of left adjoints, $\<F_{\Psi},Y_{\bbB}\>\circ\<F_{\Phi},Y_{\bbA}\>$ is left adjoint. And since $\<F_{\Phi},Y_{\bbA}\>\circ Y_{\bbA}=F_{\Phi}$, it suffices to show that  $\<F_{\Psi},Y_{\bbB}\>\circ F_{\Phi}\cong F_{\Psi\tensor\Phi}$. 
But for any $a\in\bbA$ one calculates that
\begin{eqnarray*}
\Big(\<F_{\Psi},Y_{\bbB}\>\circ F_{\Phi}\Big)(a)
 & = & \<F_{\Psi},Y_{\bbB}\>\Big(\Phi(-,a)\Big) \\
 & = & \colim\Big(\Phi(-,a),F_{\Psi}\Big) \\
 & = & \Big(\Psi\tensor_{\bbB}\Phi\Big)(-,a) \\
 & = & F_{\Psi\tensor_{\bbB}\Phi}(a).
\end{eqnarray*}
That this functor is really a 2-functor, and that it is locally an equivalence, follows from the preceding remarks, in particular (\ref{1008}).
\endofproof
Restricting to skeletal cocomplete $\Q$-categories, $\Cocont\skel(\Q)$ is a quantaloid and the 2-functor above corestricts to a quantaloid homomorphism $\Dist(\Q)\to\Cocont\skel(\Q)$ that induces sup-isomorphisms between the hom-sup-lattices.
\begin{example}
For a $\2$-category $\bbA$, i.e.~an ordered set, $\P\bbA$ is set of down-sets of $\bbA$, ordered by inclusion. The Yoneda embedding $Y_{\bbA}\: \bbA\to\P\bbA$ sends an element to the principal down-set that it determines. An eventual left adjoint to $Y_{\bbA}$ is ``taking the supremum''. The isomorphism $\Dist(\2)(\bbA,\bbB)\cong\Cocont\skel(\2)(\P\bbA,\P\bbB)$ says that an ideal relation between given orders is the same thing as a sup-morphism between the sup-lattices of down-sets of those orders. In fact, $\Cocont\skel(\2)=\Sup$. 
\end{example}

\subsection*{Covariant presheaves, free completion}

The results in the previous sections may be dualized; one obtains that ``the $\Q$-category of covariant presheaves on $\bbA$ is $\bbA$'s free completion in $\Cat(\Q)$''. Let us quickly do this exercise.
\par
A {\em covariant presheaf of type $C\in\Q$} on a $\Q$-category $\bbA$ is a distributor $\phi\:\bbA\dist*_C$. Covariant presheaves form the object-set of a $\Q$-category $\P\+\bbA$, whose hom-arrows are, for $\phi\:\bbA\dist*_C$ and $\psi\:\bbA\dist *_{C'}$,
$$\P\+\bbA(\psi,\phi)=\mbox{ single element of }\{\phi,\psi\}\mbox{ in }\Dist(\Q).
$$This $\Q$-category satisfies
\begin{equation}\label{3000}
\mbox{for every $\Q$-category $\bbC$, }\Dist(\Q)(\bbA,\bbC)\simeq\Cat(\Q)\co(\bbC,\P\+\bbA),\end{equation}
and it is complete. There is a Yoneda embedding $Y_{\bbA}\+\:\bbA\to\P\+\bbA$, sending an object $a\in\bbA$ to the representable distributor 
$\bbA(a,-)\:\bbA\dist*_{ta}$, which is fully faithful. $\bbA$ is complete if and only if $Y_{\bbA}\+$ is a left adjoint in $\Cat(\Q)$; but in any case is $Y_{\bbA}\+$ cocontinuous.
\par
Continuous functors whose domain is a complete $\Q$-category, coincide with the right adjoints. Denoting $\Cont(\Q)$ for the sub-2-category of $\Cat(\Q)$ whose objects are complete $\Q$-categories and arrows are continuous functors, we have
$$\P\+\dashv\U\:\Cont(\Q)\adjar\Cat(\Q);$$
that is, $\P\+\bbA$ is the free completion of $\bbA$ in $\Cat(\Q)$. Furthermore, $\Cont(\Q)$ has stable local limits, and there is 2-functor, which is locally an equivalence, like so: $\Dist(\Q)\to\Cont(\Q)^{\sf coop}$ sends a distributor $\Phi\:\bbA\to\bbB$ to the right Kan extension
$(F^{\Phi},Y\+_{\bbA})\:\P\+\bbB\to\P\+\bbA$, where $F^{\Phi}\:\bbB\to\P\+\bbA$ is determined by $\Phi$ through the equivalence in (\ref{3000}).
\par
Since a $\Q$-category is complete if and only if it is cocomplete, and since between such categories continuous (cocontinuous) functors coincide with right adjoints (left adjoints), there is an isomorphism of (locally cocompletely ordered) 2-categories
$$\Cocont(\Q)\doubiso\Cont(\Q)^{\sf coop}\:\Big(F\:\bbA\to\bbB\Big)\doubar\Big(F^*\:\bbB\to\bbA\Big)$$
where $F\dashv F^*$ in $\Cat(\Q)$. Considering only skeletal objects, this restricts to an isomorphism of quantaloids $\Cocont\skel(\Q)\cong\Cont^{\sf coop}\skel(\Q)$.
\par
Assembling previous results we obtain that: for a $\Q$-category $\bbA$, $Y_{\bbA}\:\bbA\to\P\bbA$ is a right adjoint if and only if $Y_{\bbA}\+\:\bbA\to\P\+\bbA$ is a left adjoint; this is precisely the case when $\bbA$ is (co)complete.
It is interesting to take a closer look at this result.
First observe that, in any case, there is an adjunction
$$\xymatrix@=15mm{\P\bbA\ar@{}[r]|{\perp}\ar@/^4mm/[r]^{[-,\bbA]} & \P\+\bbA\ar@/^4mm/[l]^{\{-,\bbA\}} }$$
in $\Cat(\Q)$, where the functors involved are obtained by calculating $\{-,-\}$ and $[-,-]$ in $\Dist(\Q)$:
$$[-,\bbA]\:\P\bbA\to\P\+\bbA\:\Big(\phi\:*_A\dist\bbA\Big)\mapsto\Big(\big[\phi,\bbA\big]\:\bbA\dist*_A\Big),$$
$$\{-,\bbA\}\:\P\+\bbA\to\P\bbA\:\Big(\psi\:\bbA\dist *_A\Big)\mapsto\Big(\big\{\psi,\bbA\big\}\:*_A\dist\bbA\Big).$$
Now considering the diagram of $\Q$-categories and functors
$$\xymatrix@R=17mm@C=12mm{
 & \bbA\ar[dl]|{Y\+_{\bbA}}\ar[dr]|{Y_{\bbA}} \\
\P\+\bbA\ar@<1mm>[rr]^{\{-,\bbA\}}\ar@{.>}@/^6mm/[ur]^R & & \P\bbA\ar@<1mm>[ll]^{[-,\bbA]}\ar@{.>}@/_6mm/[ul]_L}$$
we have on the one hand that, if $L\dashv Y_{\bbA}$ is known, then $R:=L\circ\{-,\bbA\}$ is right adjoint to $Y_{\bbA}\+$; and dually, if $Y_{\bbA}\+\dashv R$ is known, then $L:=R\circ [-,\bbA]$ is left adjoint to $Y_{\bbA}$. Actually, $L(\phi)=\colim(\phi,1_{\bbA})$ and $R(\psi)=\lim(\psi,1_{\bbA})$ (whenever one, thus also the other, exists), so that the {\em crux} of this construction is really the statement in \ref{1004}.
\begin{example} Take $\2$-categories, then
$\Cont(\2)\simeq\Inf=\Cont\skel(\2)$. For a $\2$-category $\bbA$, $\P\+\bbA$ is the set of up-sets of $\bbA$, ordered by containment (reverse inclusion); $Y_{\bbA}\+$ sends an element of $\bbA$ to the corresponding principal up-set. An eventual right adjoint to $Y_{\bbA}\+$ is ``infimum''. It is well-known that $\Sup\cong\Inf^{\sf coop}$ in $\QUANT$ by passing from left-adjoints to right-adjoints in $\Cat(\2)$. The adjunction between $\P\bbA$ (down-sets of $\bbA$, ordered by inclusion) and $\P\+\bbA$ (up-sets of $\bbA$, ordered by containment) is constituted by the left adjoint that sends a down-set to the up-set of its upper bounds, and the right adjoint that sends an up-set to the down-set of its lower bounds.
\end{example}
\begin{example}
For a $\Q$-object $X$, the category of covariant presheaves on $*_X$ is denoted $\P\+ X$ rather than $\P\+(*_X)$ (cf.~\ref{a04}).
\end{example}

\section{Cauchy completion and Morita equivalence}\label{ghi} 

\subsection*{Cauchy complete categories}

A left adjoint distributor $\Phi\:\bbC\dist\bbA$ is often called a {\em Cauchy
distributor} (into
$\bbA)$; we will systematically denote its right adjoint as $\Phi^*$:
$$\xymatrix@=15mm{\bbA\ar@{}[r]|{\perp}\ar@<1mm>@/^3mm/[r]|{\distsign}^{\Phi}&\bbB\ar@<1mm>@/^3mm/[l]|{\distsign}^{\Phi^*}}.$$
 In particular, a {\em Cauchy
presheaf} $\phi\:*_C\dist\bbA$ is a contravariant presheaf with a right adjoint $\phi^*\:\bbA\dist *_C$, whose 
elements are $\phi^*(a)\:ta\to C$.
\par
A Cauchy distributor $\Phi\:\bbC\dist\bbA$ is said to {\em converge} (to a
functor
$F\:\bbC\to\bbA$) if (there exists a functor
$F\:\bbC\to\bbA$ such that) $\Phi=\bbB(-,F-)$---or equivalently, $\Phi^*=\bbB(F-,-)$. A Cauchy distributor $\Phi\:\bbC\dist\bbA$
need not converge, but if it does then the functor $F\:\bbC\to\bbA$ that it converges to, is
essentially unique. We are interested in those
$\Q$-categories in which ``all Cauchy distributors converge'': say that $\bbA$ is {\em Cauchy complete} if, for all $\bbC$, any Cauchy distributor
$\Phi\:\bbC\dist\bbA$ converges. In other words, letting $\Map(\Dist(\Q))(\bbC,\bbA)$ denote the full subcategory of
$\Dist(\Q)(\bbC,\bbA)$ whose objects are the Cauchy distributors (``maps''), $\bbA$ is Cauchy complete if and only if, for every
$\bbC$, the mapping $$\Cat(\Q)(\bbC,\bbA)\to\Map(\Dist(\Q))(\bbC,\bbA)\:F\mapsto\bbA(-,F-)$$ is surjective; recalling that it is always ``essentially injective'', we thus ask for it to be an equivalence. 
\par
Fortunately, to have Cauchy completeness, convergence of the Cauchy presheaves suffices!
\begin{proposition}\label{134} A
$\Q$-category $\bbA$ is Cauchy complete if and only if all Cauchy presheaves on $\bbA$ converge (i.e.~are representable).
\end{proposition}
\proof We only need to prove the ``if''. Let $\Phi\dashv\Phi^*\:\bbA\adjdist\bbC$, then for each
$c\in\bbC_0$, $\Phi(-,c)\dashv\Phi^*(c,-)\:\bbA\adjdist\bbC$, because
$\Phi\tensor_{\bbC}\bbC(-,c)=\Phi(-,c)$ and $\bbC(c,-)\tensor_{\bbC}\Phi^*=\Phi^*(c,-)$, and by composition of adjoints.
By assumption there exists for each
$c\in\bbC_0$ a functor
$Fc\:*_{tc}\to\bbA$ -- that we will henceforth identify with the single element of type
$tc$ that it picks out -- such that
$\bbA(-,Fc)=\Phi(-,c)$ or equivalently, $\bbA(Fc,-)=\Phi^*(c,-)$. But the object
mapping
$\bbC_0\to\bbA_0\:c\mapsto Fc$ is functorial: since $\bbC\leq
\Phi^*\tensor_{\bbA}\Phi$ it follows that
$$\bbC(c',c)\leq\Phi^*(c',-)\tensor_{\bbA}\Phi(-,c) =\bbA(Fc',-)\tensor_{\bbA}\bbA(-,Fc)=
\bbA(Fc',Fc).$$ Clearly $\Phi\dashv\Phi^*$ converges to the functor
$F\:\bbC\to\bbA\:c\mapsto Fc$.
\endofproof
\par 
The next proposition explains Cauchy completeness as ``admitting all Cauchy-weighted
colimits''.
\begin{proposition}\label{147} A $\Q$-category is Cauchy complete if and only if it admits
all colimits weighted by a
Cauchy distributor, if and only if it admits all colimits weighted by a Cauchy (contravariant) presheaf.
\end{proposition}
\proof Let $\bbB$ be Cauchy complete; given a Cauchy distributor $\Phi\:\bbC\dist\bbA$ and a functor $F\:\bbA\to\bbB$, the $\Phi$-weighted colimit of $F$ is the functor $K\:\bbC\to\bbB$ to which the Cauchy distributor $\bbB(-,F-)\tensor_{\bbA}\Phi$ converges:
$$\bbB(K-,-)=\Phi^*\tensor_{\bbA}\bbB(F-,-)=\Big[\Phi,\bbB(F-,-)\Big].$$
Suppose now that $\bbB$ admits all colimits weighted by a Cauchy presheaf; then every such Cauchy presheaf is representable because $\phi\:*_C\dist\bbB$ (with $\phi\dashv\phi^*$) converges to 
$\colim(\phi,1_{\bbB})$:
$$\bbB\Big(\colim(\phi,1_{\bbB})-,-\Big)=\Big[\phi,\bbB(1_{\bbB}-,-)\Big]=\phi^*\tensor_{\bbB}\bbB(-,-)=\phi^*.$$
By \ref{134} $\bbB$  is Cauchy complete. (The remaining implication is trivial.)
\endofproof
\par
`Cauchy completeness' is a self-dual notion: a $\Q$-category $\bbA$ is Cauchy complete if and only if the $\Q\op$-category $\bbA\op$ is. (This really follows directly from the definition.) Therefore we have -- by dualizing the statement of \ref{147}, using that `Cauchy completeness' is a self-dual notion --  at once the following.
\begin{proposition}\label{147.1}
A $\Q$-category $\bbA$ is Cauchy complete if and only if it admits all limits weighted by a right adjoint distributor, if and only if it admits all limits weighted by a right adjoint (covariant) presheaf.
\end{proposition}
In the rest of this section we work in terms of colimits.
\begin{example}
Every $\2$-enriched category $\bbA$ is Cauchy complete: a down-set of $\bbA$ (i.e.~a contravariant presheaf on $\bbA$) is Cauchy if and only if it is principal (i.e.~a representable presheaf).
\end{example}
\begin{example}\label{a000}
The $\Rel(\Omega)$-category associated to a sheaf on a locale $\Omega$ (see \ref{a05.0}) is (skeletal and) Cauchy complete. In fact, it is also {\em symmetric} in the sense that it equals its opposite (which makes sense because $\Rel(\Omega)\op=\Rel(\Omega)$). [Walters, 1981] proves that sheaves on $\Omega$ are {\em precisely} the symmetric, skeletal, Cauchy complete $\Rel(\Omega)$-categories. Dropping symmetry and skeletality, Cauchy complete $\Rel(\Omega)$-categories are ordered sheaves on $\Omega$ [Betti\etal, 1983; Stubbe, 2004c].
\end{example}
\begin{example}
Consider again Lawvere's generalized metric spaces (see \ref{a06}); let $\bbA$ denote such a space. In [Lawvere, 1973] it is proved that a left adjoint distributor from a one-point space into $\bbA$ is precisely (an equivalence class of) a Cauchy sequence in $\bbA$---and it is therefore that left adjoint distributors are also called Cauchy distributors. Cauchy completeness of $\bbA$ as enriched category coincides with its Cauchy completeness as metric space (``all Cauchy sequences/distributors converge'').
\end{example}

\subsection*{Cauchy completing a category}

Not every $\Q$-category $\bbA$ is Cauchy
complete; still we are interested in ``classifying'' the Cauchy distributors on $\bbA$: we want to construct a $\Q$-category $\bbA\cc$ such that
\begin{equation}\label{138}
\mbox{for every $\Q$-category $\bbC$, }\Map(\Dist(\Q))(\bbC,\bbA)\simeq\Cat(\Q)(\bbC,\bbA\cc).
\end{equation}   Much like the
construction of $\P\bbA$ we have the following proposition (the proof of which is much like
that of \ref{102}, so it is omitted).
\begin{proposition}\label{139} For any $\Q$-category
$\bbA$, the ``Cauchy completion $\bbA\cc$ of $\bbA$'' is, by definition, the full subcategory of
$\P\bbA$ determined by the Cauchy presheaves on $\bbA$:
$$(\bbA\cc)_0=\{\phi\in(\P\bbA)_0\mid
\phi\mbox{ is Cauchy}\}.$$ This $\Q$-category is skeletal and it satisfies (\ref{138}):
\begin{itemize}
\item for a Cauchy distributor $\Phi\:\bbC\dist\bbA$ the corresponding functor $F_{\Phi}\:\bbC\to\bbA\cc$ maps $c\in\bbC_0$ onto the Cauchy presheaf $\Phi(-,c)\:tc\dist\bbA$;
\item and for a functor $F\:\bbC\to\bbA\cc$ the corresponding Cauchy distributor $\Phi_F\:\bbC\dist\bbA$ has elements, for $(c,a)\in\bbC_0\times\bbA_0$, $\Phi_F(a,c)=F(c)(a)$ (and its right adjoint has elements $\Phi_F^*(c,a)=F(c)^*(a)$).
\end{itemize}
\end{proposition}
Note that the hom-arrows in $\bbA\cc$ are identical to those in
$\P\bbA$, but since objects
$\phi,\psi\in\bbA\cc$ have right adjoints in $\Dist(\Q)$ we may write (with slight abuse
of notation) that $\bbA\cc(\psi,\phi)= \psi^*\tensor_{\bbA}\phi$.
This will be a helpful trick in many of the calculations that follow.
\par
Plugging the identity distributor $\bbA\:\bbA\dist\bbA$ into the equivalence (\ref{138}) gives us a functor
$$i_{\bbA}\:\bbA\to\bbA\cc\:a\mapsto\Big(\bbA(-,a)\:*_{ta}\dist\bbA\Big)$$ 
which is precisely the factorization of the Yoneda embedding $Y_{\bbA}\:\bbA\to\P\bbA$ over the full
inclusion
$j_{\bbA}\:\bbA\cc\to\P\bbA$ that defined $\bbA\cc$. An alternative for the classifying property in (\ref{138}) is then: for any $\bbC$, a distributor
$\Phi\:\bbC\dist\bbA$ is Cauchy
if and only if the corresponding functor
$F_{\Phi}\:\bbC\to\P\bbA$ (as in~\ref{102}) factors through
$j_{\bbA}\:\bbA\cc\to\P\bbA$.
It follows too that $i_{\bbA}$ enjoys properties similar to those of $Y_{\bbA}$ (especially \ref{104} and \ref{109}).
\begin{proposition}\label{139.1}
For any Cauchy presheaf $\phi\:*_C\dist\bbA$ and any object $a\in\bbA_0$, $\bbA\cc(i_{\bbA}(a),\phi)=\phi(a)$ and $\bbA\cc(\phi,i_{\bbA}(a))=\phi^*(a)$. In particular is $i_{\bbA}\:\bbA\to\bbA\cc$ fully faithful.
\end{proposition}
\begin{proposition}\label{139.2} Every Cauchy presheaf
on a $\Q$-category $\bbA$ is canonically the colimit of representables in $\bbA\cc$: for $\phi\:*_C\dist\bbA$ Cauchy,
the $\phi$-weighted colimit of $i_{\bbA}\:\bbA\to\bbA\cc$ (exists and) is $\phi$ itself. In other words, $\colim(\phi,i_{\bbA})=\phi$ (where we may write an equality instead of merely an isomorphism, for $\bbA\cc$ is skeletal).
\end{proposition}
\par
Up to now we have encountered many similarities between $Y_{\bbA}\:\bbA\to\P\bbA$ and $i_{\bbA}\:\bbA\to\bbA\cc$; but here are two remarkable differences.
\begin{proposition}\label{150.0}
A $\Q$-category $\bbA$ is Cauchy complete if and only if $i_{\bbA}\:\bbA\to\bbA\cc$ is an essentially surjective functor, if and only if it is an equivalence in $\Cat(\Q)$\footnote{Cocompleteness of $\bbA$ is in general not equivalent to $Y_{\bbA}\:\bbA\to\P\bbA$ being an equivalence in $\Cat(\Q)$!}.
\end{proposition}
\proof Since $i_{\bbA}\:\bbA\to\bbA\cc$ is fully faithful, it is clear that it is an equivalence if and only if it is essentially surjective. But the essential surjectivity, saying that any $\phi\in\bbA\cc$ is isomorphic to $i_{\bbA}(a)$ for some $a\in\bbA$, is precisely the statement in \ref{134}.
\endofproof
\begin{proposition}\label{150} Any $\Q$-category $\bbA$ is isomorphic to its Cauchy
completion in the quantaloid $\Dist(\Q)$\footnote{But $\bbA$ needn't be isomorphic to $\P\bbA$ in $\Dist(\Q)$!}. 
\end{proposition}
\proof $i_{\bbA}\:\bbA\to\bbA\cc$ is a fully faithful
functor, so the unit of the induced adjoint pair of distributors
$\bbA\cc(-,i_{\bbA}-)\dashv\bbA\cc(i_{\bbA}-,-)$ is an equality. It suffices to prove that also its
co-unit is an equality. Thereto, consider $\phi,\psi\in(\bbA\cc)_0$
and calculate that 
\begin{eqnarray*}
\bbA\cc(\psi,i_{\bbA}-)\tensor_{\bbA}\bbA\cc(i_{\bbA}-,\phi)
 & = & \bigvee_{a\in\bbA_0}\bbA\cc\Big(\psi,\bbA(-,a)\Big)\circ\bbA\cc\Big(\bbA(-,a),\phi\Big) \\
 & = &
\psi^*\tensor_{\bbA}\bigvee_{a\in\bbA_0}\Big{(}\bbA(-,a)\circ\bbA(a,-)\Big{)\tensor_{\bbA}\phi}
\\
 & = & \psi^*\tensor_{\bbA}\bbA\tensor_{\bbA}\phi \\
 & = & \psi^*\tensor_{\bbA}\phi \\
 & = & \bbA\cc(\psi,\phi).
\end{eqnarray*} That is to say,
$\bbA\cc(-,i_{\bbA}-)\tensor_{\bbA}\bbA\cc(i_{\bbA}-,-)=\bbA\cc$, as wanted.
\endofproof
A direct consequence of the above is that $\Dist(\Q)$ is equivalent to its full subquantaloid $\Dist\cc(\Q)$ determined by the Cauchy complete categories.
\par
The following result, although it uses \ref{150} in its proof, may be considered analogous to \ref{107}. 
\begin{proposition}\label{146} For any $\bbA$, the Cauchy completion
$\bbA\cc$ is Cauchy complete.
\end{proposition} 
\proof
It suffices to prove that any Cauchy presheaf on $\bbA\cc$ is representable: for any Cauchy $\Phi\:*_C\dist\bbA\cc$ we must find a Cauchy $\phi\:*_C\dist\bbA$ such that $\bbA\cc(-,\phi)=\Phi$. Using $\bbA\cong\bbA\cc$ in $\Dist(\Q)$ (see \ref{150}) we put $\phi:=\bbA\cc(i_{\bbA}-,-)\tensor_{\bbA\cc}\Phi$, $$\xymatrix@=15mm{\ast_C\ar[r]|{\distsign}^{\Phi}\ar@/_2pc/[rr]|{\distsign}_{\phi} & \bbA\cc\ar[r]|{\distsign}^{\bbA\cc(i_{\bbA}-,-)} & \bbA,}$$
which is Cauchy because it is the composition of Cauchy distributors ($\bbA\cc(i_{\bbA}-,-)$ is not only right but also left adjoint because it is invertible). Now for any $\psi\in\bbB\cc$,
\begin{eqnarray*}
\bbA\cc(\psi,\phi)
& = & \psi^*\tensor_{\bbA}\phi \\
& = & \psi^*\tensor_{\bbA}\bbA\cc(i_{\bbA}-,-)\tensor_{\bbA\cc}\Phi \\
& = & \bbA\cc(\psi,i_{\bbA}-)\tensor_{\bbA}\bbA\cc(i_{\bbA}-,-)\tensor_{\bbA\cc}\Phi \\
& = & \bbA\cc(\psi,-)\tensor_{\bbA\cc}\Phi \\
& = & \Phi(\psi),
\end{eqnarray*}
so $\phi$ indeed represents $\Phi$.
\endofproof
In view of the examples, in particular \ref{a000}, it is tempting to think of the Cauchy completion $\bbA\mapsto \bbA\cc$ of a $\Q$-category as a sheafification. The following will only strengthen that intuition.

\subsection*{The universality of the Cauchy completion}

The proof of the
following is analogous to that of \ref{122}, and therefore it is omitted.
\begin{proposition}\label{148.0}
A functor $F\:\bbA\to\bbB$ has a pointwise left Kan extension along $i_{\bbA}\:\bbA\to\bbA\cc$ if and only if, for every $\phi\in\bbA\cc$, the $\phi$-weighted colimit of $F$ exists; in this case $\<F,i_{\bbA}\>\circ i_{\bbA}\cong F$\footnote{But it is not true in general that $\<F,i_{\bbA}\>\:\bbA\cc\to\bbB$ is a left adjoint!}. Actually, $\<F,i_{\bbA}\>$ is the essentially unique factorization of $F$ through $i_{\bbA}$.
\end{proposition}
In other words, denoting $\Cat\cc(\Q)$ for the full sub-2-category of $\Cat(\Q)$ determined by the Cauchy complete categories, we have for $\Q$-categories $\bbA$ and $\bbB$, with $\bbB$ Cauchy
complete, an equivalence induced by composition with
$i_{\bbA}\:\bbA\to\bbA\cc$
$$\Cat\cc(\Q)(\bbA\cc,\bbB)\iso
\Cat(\Q)(\bbA,\bbB) \: G \mapsto G\circ i_{\bbA}$$ (with inverse $F\mapsto\<F,i_{\bbA}\>$).
\begin{proposition}\label{148}
Cauchy completing provides a left (bi)adjoint to the embedding of $\Cat\cc(\Q)$ in $\Cat(\Q)$: $(-)\cc\dashv\V\:\Cat\cc(\Q)\adjar\Cat(\Q)$.
\end{proposition}
The 2-functor
$$\Cat\cc(\Q)\to\Map(\Dist\cc(\Q))\:\Big(F\:\bbA\to\bbB\Big)\mapsto\Big(\bbB(-,F-)\:\bbA\dist\bbB\Big)$$
is the identity on objects, and induces equivalences of the hom-orders: it is a biequivalence.
Since $\Dist(\Q)$ itself is equivalent to $\Dist\cc(\Q)$, we may record that
$$\Cat\cc(\Q)\simeq\Map(\Dist\cc(\Q))\simeq\Map(\Dist(\Q)).$$
\par Finally we wish to discuss the ``interplay'' between $\bbA$, $\P\bbA$, $\P\+\bbA$ and $\bbA\cc$.
\begin{proposition}\label{155}
For $\Q$-categories $\bbA$ and
$\bbB$, the following are equivalent: 
\begin{enumerate}
\item $\bbA\cong\bbB$ in $\Dist(\Q)$;
\item $\bbA\cc\simeq\bbB\cc$ in $\Cat(\Q)$;
\item $\P\bbA\simeq\P\bbB$ in $\Cat(\Q)$;
\item $\P\+\bbA\simeq\P\+\bbB$ in $\Cat(\Q)$.
\end{enumerate}
\end{proposition}
\proof
$(1\Leftrightarrow 2)$ $\bbA\cong\bbB$ in $\Dist(\Q)$ implies and is implied by $\bbA\cc\cong\bbB\cc$ in $\Dist\cc(\Q)$ (using \ref{139}), which by the biequivalence $\Cat\cc(\Q)\simeq\Map(\Dist\cc(\Q))$ means that $\bbA\cc\simeq\bbB\cc$ in $\Cat\cc(\Q)$, or in $\Cat(\Q)$.
\par
$(1\Leftrightarrow 3)$ $\bbA\cong\bbB$ in $\Dist(\Q)$ if and only if $\P\bbA\simeq\P\bbB$ in $\Cocont(\Q)$ by the fully faithful, locally fully faithful 2-functor $\Dist(\Q)\to\Cocont(\Q)$ of \ref{132}; but to give an equivalence between cocomplete categories in $\Cocont(\Q)$ is the same as in $\Cat(\Q)$ (because the functors constituting the equivalence are always left adjoint, thus cocontinuous).
\par
$(1\Leftrightarrow 4)$ is dual to the previous argument.
\endofproof
Whenever for two
$\Q$-categories $\bbA$ and
$\bbB$ the equivalent conditions of \ref{155} are fulfilled, then
$\bbA$ is {\em Morita equivalent} to
$\bbB$.
The morality is then that ``a $\Q$-category is but a
presentation of its Cauchy completion''; in particular, if one only cares about
$\Q$-category theory ``up to Cauchy completion'' (i.e.~``up to Morita equivalence''), then
one can forget about functors and Cauchy completion and Morita equivalence altogether, and
treat $\Map(\Dist(\Q))$ as ``the 2-category of $\Q$-categories and functors''.

\section{Appendix: Distributor calculus}

\subsection*{Lax (co)limits in a quantaloid}

Let $\Q$ denote a quantaloid, and $\D$ a small category. A {\em
lax functor}
$F\:\D\to\Q$ is a mapping of objects
$\D_0\to\Q_0\:D\mapsto FD$ together with an action on hom-sets
$\D(D,D')\to\Q(FD,FD')\:f\mapsto Ff$ such that 
\begin{itemize}
\item for all $D\in\D_0$, $1_{FD}\leq F1_D$;
\item for all composable $f,g\in\D$, $(Fg\circ Ff)\leq F(g\circ f)$.
\end{itemize}  Given two such lax functors $F,G\:\D\biar\Q$, a family 
$(\theta_D\:FD\to GD)_{D\in\D_0}$ of
$\Q$-arrows constitutes a {\em lax-natural transformation}
$\theta\:F\tto G$ if
\begin{itemize}
\item  for all $f\:D\to D'$ in $\D$, $(\theta_{D'}\circ Ff)\geq(Gf\circ\theta_D)$;
\end{itemize} and the family is said to constitute an {\em oplax-natural transformation} if
\begin{itemize}
\item  for all $f\:D\to D'$ in $\D$, $(\theta_{D'}\circ Ff)\leq(Gf\circ\theta_D)$.
\end{itemize}
The collection of lax transfos, resp.~oplax transfos, between two parallel
lax functors
$F,G\:\D\biar\Q$ constitutes a sup-lattice (with ``componentwise supremum''); we'll denote
it as
$\Lax(F,G)$, resp.~$\OpLax(F,G)$. Lax transfos (oplax transfos) may also be composed: this
too is done ``componentwise''. The notions of `lax transfo' and `oplax transfo' are dual in
the following sense: two lax functors
$F,G\:\D\biar\Q$ may be thought of as
$F\op,G\op\:\D\op\biar\Q\op$, and then we have that
$\Lax(F,G)=\OpLax(G\op,F\op)$.
\par Any (ordinary) functor $F\:\D\to\Q$ is trivially lax, and any (ordinary) natural
transformation
$\alpha\: F\tto G$ between functors is trivially (op)lax. The set of natural transformations
between two functors $F,G\:\D\biar\Q$ is a sub-sup-lattice of both the sup-lattices
$\Lax(F,G)$ and $\OpLax(F,G)$ (it is their intersection).
\par In what follows, $\Delta_A\:\D\to\Q$ denotes the constant functor at the object
$A\in\Q$ (on whatever category $\D$). For an arrow $k\:A\to B$ in $\Q$ the constant
transformation from
$\Delta_A$ to
$\Delta_B$ will be denoted as $\Delta_k$.
\begin{definition}\label{9} Let $\D$ be a small category, and $\Q$ any quantaloid. For a lax
functor
$F\:\D\to\Q$ a {\em lax cone} over $F$ consists of an object $L\in\Q$ and a lax transfo
$\pi\in\Lax(\Delta_L,F)$. Such a lax cone $(L,\pi)$ is the (necessarily essentially unique)
{\em lax conical limit} of $F$ if it is universal: for all
$A\in\Q$ there is an isomorphism\footnote{Actually the definition of ``lax (co)limit'' in a
general bicategory asks only for an {\em equivalence}, and not for an isomorphism. But an
equivalence between antisymmetric orders is always an isomorphism, so we put that right
away in this definition.} of sup-lattices induced by composition with
$\pi\:
\Delta_L\tto F$, 
\begin{equation}\label{10.1}
\Q(A,L)\iso\Lax(\Delta_A,F)\:k\mapsto\pi\circ\Delta_k.
\end{equation} {\em Lax cocone} and {\em lax conical colimit} are the dual notions: a lax cocone on
$F$ is an oplax transfo
$\sigma\in\OpLax(F,\Delta_C)$, and it is the lax colimit if it is universal.
\end{definition}  So a lax functor $F\:\D\to\Q$ admits a lax limit (lax colimit) in $\Q$ if
and only if
$F\op\:\D\op\to\Q\op$ admits a lax colimit (lax limit) in $\Q\op$---this is the duality
between limits and colimits.
\par In general, for an ordinary functor
$F\:\D\to\Q$, the limit on $F$ is not necessarily its lax limit: whereas any cone
$\alpha\in\Nat(\Delta_L,F)$ is trivially also a lax cone, it is not true that the universal
property for such an $\alpha$ to be a {\em limit} implies the universal property for
$\alpha$ to be a {\em lax} limit. But in some cases lax limits do coincide with ordinary
ones. (And by duality we can apply the result to lax colimits, of course.)
\begin{proposition}\label{19} Suppose that, for some lax functor
$F_1\:\D_1\to\Q$ and some (ordinary) functor
$F_2\:\D_2\to\Q$ into a quantaloid $\Q$, for all $X\in\Q$ there is a sup-isomorphism
\begin{equation}\label{20}
\Lax(\Delta_X,F_1)\iso \Nat(\Delta_X,F_2)
\end{equation} which is natural in $X$. Then some $\pi\in\Lax(\Delta_L,F_1)$ is the lax
limit of $F_1$ if and only if $\pi'\in\Nat(\Delta_L,F_2)$, the image of $\pi$ under the
isomorphism above (taking $X=L$), is the limit of
$F_2$.
\end{proposition}
\proof
To say that the functor $F_2\:\D_2\to\Q$ has limit $\pi'\in\Nat(\Delta_L,F_2)$ is to say that there are bijections
$$\beta'_X\:\Q(X,L)\iso\Nat(\Delta_X,F_2)\:f\mapsto \pi'\circ\Delta_f$$
defining an isomorphism of the $\Set$-valued functors $\Q(-,L)$ and $\Nat(\Delta_{-},F_2)$ on (the underlying category of) $\Q\op$ in $\CAT$; and $\pi'=\beta'_L(1_L)$. However, not only is each $\Nat(\Delta_X,F_2)$ is a sup-lattice, $\Nat(\Delta_{-},F_2)$ also preserves suprema of $\Q$-arrows; it is thus a $\Sup$-valued homomorphism on $\Q\op$. Moreover, the natural bijections above -- since they are induced by composition in $\Q$ -- always preserve suprema. Therefore, $F_2$ has limit $\pi'$ if and only if $\beta'\:\Q(-,L)\tto\Nat(\Delta_{-},F_2)$ is an isomorphism of $\Sup$-valued homomorphisms on $\Q\op$ in $\QUANT$.
\par
On the other hand, to say that the lax functor $F_1\:\D_1\to\Q$ has lax limit $\pi\in\Lax(\Delta_X,F_1)$, is to say that there are sup-isomorphisms
$$\beta_X\:\Q(X,L)\iso\Lax(\Delta_X,F_1)\:f\mapsto \pi\circ\Delta_f$$
defining an isomorphism $\beta\:\Q(-,L)\tto\Lax(\Delta_{-},F_1)$ of $\Sup$-valued homomorphisms on $\Q\op$ in $\QUANT$; and $\pi=\beta_L(1_L)$.
\par
The result above then follows by composition of isomorphisms of $\Sup$-valued homomorphisms on $\Q\op$ in $\QUANT$.
\endofproof

\subsection*{Two important examples}

A set-indexed family
$(A_i)_{i\in I}$ of objects in a quantaloid $\Q$ can be identified with (the image of) the
functor $F\:I\to
\Q:i\mapsto A_i$
 whose domain is the index-set regarded as discrete category.
\begin{proposition}\label{25} For a set-indexed family
$(A_i)_{i\in I}$ of objects in a quantaloid $\Q$, the following are equivalent:
\begin{enumerate}
\item $(A_i)_{i\in I}$ has a limit (i.e.~product) in $\Q$;
\item $(A_i)_{i\in I}$ has a colimit (i.e.~coproduct) in $\Q$;
\item $(A_i)_{i\in I}$ has a lax limit in $\Q$;
\item $(A_i)_{i\in I}$ has a lax colimit in $\Q$;
\item there exist an object $A$ and arrows 
$$(\xymatrix{A\ar[r]^{p_i} & A_i\ar[r]^{s_i} & A})_{i\in I}$$ in
$\Q$ such that, for all $i,j\in I$, $p_j\circ s_i=\delta_{A_i,A_j}$ and $\bigvee_{i\in I}(s_i\circ p_i)=1_A$.
\end{enumerate} The cone $(p_i\:A\to A_i)_{i\in I}$ is then the (lax) limit of the given family, and the cocone $(s_i\:A_i\to A)_{i\in I}$
is its (lax) colimit. For each $i\in I$ the coprojection
$s_i\:A_i\to A$ is left adjoint to the projection $p_i\:A\to A_i$.
\end{proposition}
\proof The equivalence $(1\Leftrightarrow 2\Leftrightarrow 5)$ is well-known. 
$(1\Leftrightarrow 3)$ and
$(2 \Leftrightarrow 4)$ follow from 
\ref{19}: because the
diagram $F\:I\to\Q\:i\mapsto A_i$ is discrete, any lax cone on $F$ is in fact an ordinary cone; and any lax cocone is an ordinary cocone. That is to say, in this particular case we have for any
object
$A\in\Q$ that $\Nat(\Delta_A,F)=\Lax(\Delta_A,F)$ and
$\Nat(F,\Delta_A)=\OpLax(F,\Delta_A)$. The adjunctions $s_i\dashv p_i$ hold since
$1_{A_i}=p_i\circ s_i$ and $s_i\circ p_i\leq 1_A$.
\endofproof  As usual, we will refer to limits-colimits (in a quantaloid
$\Q$) as {\em direct sums} and we write them as
$\oplus_iA_i$.
\par
A monad in a quantaloid
$\Q$ is an endo-arrow
$A\Endoar{t}$ such that
$t\circ t\leq t$ and
$1_A\leq t$. It can be identified with (the image of) a lax functor
$F\:{\sf Ob}\to \Q\:*\mapsto A, 1_*\mapsto t$. Here ${\sf Ob}$ stands for the ``generic
object'': the category with one object
$*$ and one arrow $1_*$. On the other hand, $A\Endoar{t}$ is a diagram in $\Q$, thus (the image of) an ordinary
functor $F'\:{\sf Ar}\to\Q:*_1\mapsto A, *_2\mapsto A,
\alpha\mapsto t$. Now
${\sf Ar}$ stands for the ``generic arrow'': the category with two objects,
$*_1$ and $*_2$, and exactly one arrow that is not an identity,
$\alpha\:*_1\to *_2$.
\begin{proposition}\label{38} For a monad $A\Endoar{t}$ in a quantaloid
$\Q$, the following are equivalent:
\begin{enumerate}
\item $A\Endoar{t}$ has a limit (i.e.~$t,1_A\:A\biar A$ have an equalizer) in
$\Q$;
\item $A\Endoar{t}$ has a colimit (i.e.~$t,1_A\:A\biar A$ have a coequalizer) in
$\Q$;
\item $A\Endoar{t}$ has a lax limit in $\Q$;
\item $A\Endoar{t}$ has a lax colimit in $\Q$;
\item there exist an object $B$ and arrows
$$\xymatrix{B\ar[r]^p & A\ar[r]^{s}&B}$$ in $\Q$ such that $p\circ s=t$ and $s\circ p = 1_B$.
\end{enumerate} If the above is the case for the given monad $t$, then
$p\:B\to A$ is actually its (lax) limit, and $s\:A\to B$ its (lax) colimit. And the
coprojection
$s\:A\to B$ is left adjoint to the projection
$p\:B\to A$.
\end{proposition}
\proof
$(1\Leftrightarrow 2\Leftrightarrow 5)$ is well-known for idempotents in a category, so in particular applies to
monads in a quantaloid (it is immediate that a monad in $\Q$ is an idempotent).
$(1\Leftrightarrow 3)$ and $(2\Leftrightarrow 4)$ follow from~\ref{19}:
because $1_A\leq t$, it is the same to give a lax cone on
$F\:{\sf Ob}\to\Q$ as to give an ordinary cone on $F'\:{\sf Ar}\to\Q$; and lax cocones on
$F$ coincide with cocones on $F'$. More precisely, we have for any object
$X\in\Q$ identical sup-lattices
$\Nat(\Delta_X,F')=\Lax(\Delta_X,F)$ and $\Nat(F',\Delta_X)=\OpLax(F,\Delta_X)$.
The adjunction $s\dashv p$ follows because $1_A\leq t=p\circ s$ and $s\circ p=1_B$.
\endofproof 
The lax limit of a monad $A\Endoar{t}$ in a 2-category is called the {\em object of
$t$-algebras}, or the {\em Eilenberg-Moore object} for $t$, and often written as $A^t$. And the lax colimit
is the {\em object of free
$t$-algebras}, or the {\em Kleisli object}, written $A_t$. In a
quantaloid
$\Q$, $A^t$ is isomorphic to $A_t$ as soon as one (and thus also the other) exists: ``every
$t$-algebra is free''! When $\Q$ admits (free) algebra objects for all of its monads, we will say that {\em all monads in $\Q$ split}.

\subsection*{Universal property of the distributor calculus}

For the following we refer also to [Betti\etal, 1983; Carboni\etal, 1987].
\par
Given a quantaloid $\Q$, there is a quantaloid
$\Matr(\Q)$ of ``matrices with elements in $\Q$'', with:
\begin{itemize}
\item objects: $\Q_0$-typed sets\footnote{Recall that such a $\Q_0$-typed set $(X,t)$ is really just a way of writing a small discrete diagram $(tx)_{x\in X}$ in $\Q$.};
\item arrows: for two objects $(X,t)$,
$(Y,t)$, an arrow $\bbM\: (X,t)\to(Y,t)$ is a {\em matrix} of $\Q$-arrows
$\bbM(y,x)\: tx\to ty$, one for each $(x,y)\in X\times Y$;
\item two-cells:  for $\bbM,
\bbN\:(X,t)\biar(Y,t)$, $\bbM\leq\bbN$ when
$\bbM(y,x)\leq\bbN(y,x)$ for all $x$, $y$;
\item composition: for
$\bbM\:(X,t)\to(Y,t)$ and $\bbN\:(Y,t)\to(Z,t)$, the elements of $\bbN\circ\bbM\:(X,t)\to(Y,t)$ are
$(\bbN\circ\bbM)(z,x)=\bigvee_{y\in Y}\bbN(z,y)\circ\bbM(y,x)$;
\item identities: the identity matrix $\bbI_{(X,t)}\: (X,t)\to(Y,t)$ has elements $\bbI_{(X,t)}(x',x)=\delta_{tx,tx'}$ (``Kronecker deltas'').
\end{itemize}
Clearly, $\Q$ can be embedded in $\Matr(\Q)$: every arrow
$f\:A\to B$ in
$\Q$ may be viewed as a one-element matrix. We write this embedding as
$i_{\Q}\:\Q\to\Matr(\Q)$. Obviously, $i_{\Q}\:\Q\to\Matr(\Q)$ is a fully faithful homomorphism of quantaloids; it is an
equivalence of quantaloids if and only if $\Q$ has all direct sums. Actually, $\Matr(\Q)$ itself has all direct sums (so ``taking matrices'' is an essentially idempotent process), and $i_{\Q}\:\Q\to\Matr(\Q)$ describes the universal direct-sum completion of $\Q$ in
$\QUANT$.
\par
Given a quantaloid $\Q$, the following defines a quantaloid
$\Bim(\Q)$ of ``bimodules in
$\Q$'':
\begin{itemize}
\item objects: monads in
$\Q$;
\item arrows: an arrow $b\: t \dist s$ between monads $A\Endoar{t}$ and $B\Endoar{s}$ is a
{\em bimodule} from
$t$ to $s$, i.e.~a $\Q$-arrow
$b\: A\to B$ such that $s\circ b\leq b$ and $b\circ t\leq b$;
\item two-cells: put $b\leq c$ in $\Bim(\Q)(t,s)$ whenever $b\leq c$ in
$\Q(A, B)$;
\item composition: the composition $c\tensor_{s}b\:t\dist r$ of two bimodules $b\:t\dist s$
and $c\:s\dist r$ is $c\tensor_s b=c\circ b$ in
$\Q$;
\item identity: on a monad $A\Endoar{t}$ the identity bimodule is $t\: t\dist t$ itself.
\end{itemize} 
The ``tensor'' notation for composition of bimodules helps to remind that the identity on a monad $A\Endoar{t}$ is not $1_A$ but rather $t$ itself.
\par An arrow
$f\:A\to B$ may be regarded as bimodule from
$1_A$ to $1_B$; denote the corresponding embedding of quantaloids as
$j_{\Q}\:\Q\to\Bim(\Q)$. Then $j_{\Q}\:\Q\to\Bim(\Q)$ is a fully faithful homomorphism of quantaloids which is an
equivalence of quantaloids if and only if in $\Q$ all monads split. In the quantaloid $\Bim(\Q)$ all monads split (so ``taking monads and bimodules'' is an essentially idempotent process), actually, $j_{\Q}\:\Q\to\Bim(\Q)$ describes the universal split-monad completion of
$\Q$ in $\QUANT$.
\par
It is now a matter of fact that, for any quantaloid $\Q$,
$\Dist(\Q)=\Bim(\Matr(\Q))$. Indeed, take for example the definition of ``$\Q$-category'': such an $\bbA$ consists of an endo-matrix of hom-arrows $\bbA(a',a)\:ta\to ta'$ on the $\Q_0$-typed set $(\bbA_0,t)$ of objects; these data ought to satisfy, for all $a,a',a''\in\bbA_0$,
$$\bbA(a'',a')\circ\bbA(a',a)\leq\bbA(a'',a)\mbox{ and } 1_{ta}\leq\bbA(a,a),$$
or equivalently, 
$$\bigvee_{a'\in\bbA_0}\bbA(a'',a')\circ\bbA(a',a)\leq\bbA(a'',a)\mbox{ and } \delta_{ta,ta'}\leq\bbA(a',a),$$
or still equivalently, but now in terms of matrices,
$\bbA\circ\bbA\leq\bbA$ and $\bbI_{\bbA_0}\leq\bbA$.
This says precisely that $\bbA$ is a monad in $\Matr(\Q)$, thus an object of $\Bim(\Matr(\Q))$. (Similar for distributors.)
\par
Consequently we now know the
universal property of $\Dist(\Q)$: the embedding 
$k_{\Q}\:\Q\to\Dist(\Q)$ which equals $j_{\Matr(\Q)}\circ i_{\Q}$ is a fully faithful
homomorphism of quantaloids which is an equivalence if and only if $\Q$ has all direct sums
and all monads split. $\Dist(\Q)$ itself has all direct sums and all monads split (``taking categories and distributors'' is an essentially idempotent process); actually, $k_{\Q}\:\Q\to\Dist(\Q)$ describes the universal direct-sum-and-split-monad
completion of $\Q$ in $\QUANT$.
\par
Having all direct sums and a splitting for all monads, is enough for a quantaloid $\Q$ to admit all lax (co)limits.
\begin{proposition}\label{40.0}
For any quantaloid $\Q$, the following are equivalent:
\begin{enumerate}
\item $\Q$ has direct sums and all monads split;
\item $\Q$ has all lax limits;
\item $\Q$ has all lax colimits;
\item for any lax functor $F\:\D\to\Q$ there exist an object $L$ and arrows $$(\xymatrix{L\ar[r]^{p_D}&FD\ar[r]^{s_D}&L})_{D\in\D_0}$$
in $\Q$ such that 
$$\bigvee\{s_D\circ p_D\mid D\in\D_0\}=1_L\mbox{ and }p_{D'}\circ s_D=\bigvee\{Fd\mid d\in\D(D,D')\}.$$
\end{enumerate}
With notations as in the last sentence, the $(p_D\:L\to FD)_{D\in\D_0}$ form the lax limit of $F$, and the $(s_D\:FD\to L)_{D\in\D_0}$ its lax colimit; moreover, for every $D\in\D_0$, the coprojection $s_D$ is left adjoint to the projection $p_D$.
\end{proposition}
\proof
$(2\Rightarrow 1)$ and $(3\Rightarrow 1)$ follow from \ref{25} and \ref{38}.
\par
For $(4\Rightarrow 2)$ and $(4\Rightarrow 3)$, first observe that $s_D\dashv p_D$: on the one hand $1_L=\bigvee\{s_D\circ p_D\mid D\in\D_0\}\geq s_D\circ p_D$, and on the other $p_D\circ s_D=\bigvee\{Fd\mid d\in\D(D,D')\}\geq F(1_D)\geq 1_{FD}$. Knowing this, it follows from $p_{D'}\circ s_D=\bigvee\{Fd\mid d\in\D(D,D')\}$ that for every $d\:D\to D'$, $p_{D'}\circ s_D\geq Fd$, hence also $s_D\geq s_{D'}\circ Fd$ and $p_{D'}\geq Fd\circ p_D$. Hence $\pi=(p_D)_D$ is a lax cone, and $\sigma=(s_D)_D$ a lax cocone, on $F$. As for their universality, we must verify that for any $X\in\Q$ the sup-morphisms
$$\Q(X,L)\to\Lax(\Delta_X,F)\:x\mapsto \pi\circ\Delta_x,$$
$$\Q(L,X)\to\Lax(F,\Delta_X)\:x\mapsto\Delta_x\circ\sigma$$
are bijective (in which case they are sup-isomorphisms). But a calculation shows that their respective inverses are
$$\pi'=(p'_D)_D\mapsto \bigvee\{s_D\circ p'_D\mid D\in\D_0\},$$
$$\sigma'=(s'_D)_D\mapsto \bigvee\{s'_D\circ p_D\mid D\in\D_0\}.$$
\par
For $(1\Rightarrow 4)$, note that a lax functor $F\:\D\to\Q$ determines a $\Q$-category $\bbD$ with objects $\bbD_0=\D_0$, types $tD=FD$, and hom-arrows $\bbD(D',D)=\bigvee\{Fd\mid d\in\D(D,D')\}$: the inequalities $\bbD(D'',D')\circ\bbD(D',D)\leq\bbD(D'',D)$ and $1_{FD}\leq\bbD(D,D)$ hold precisely because $F$ is lax! Since $\Q$ has all direct sums and all monads split, this $\Q$-category $\bbD$ must be isomorphic in $\Dist(\Q)$ to a one-object $\Q$-category, call it $*_L$. This means that there exist distributors $\sigma\:\bbD\dist*_L$ and $\pi\:*_L\dist\bbD$ such that $\sigma\tensor_{*_L}\pi=\bbD$ and $\pi\tensor_{\bbD}\sigma=1_L$. 
Such $\sigma$ and $\pi$ are really collections of $\Q$-arrows 
$$\Big(\xymatrix@C=20mm{L\ar[r]^{p_D:=\pi(D)} & FD\ar[r]^{s_D:=\sigma(D)} & L}\Big)_{D\in\D_0}$$
satisfying 
$$p_{D'}\circ s_D=\bbD(D',D)=\bigvee\{Fd\mid d\in\D(D,D')\},$$
$$\bigvee\{s_D\circ p_D\mid D\in\bbD_0\}=\bigvee\{s_D\circ p_D\mid D\in\D_0\}=1_L.$$
\endofproof

\end{document}